\documentclass[11pt]{article}

\newcommand{\hide}[1]{}

\usepackage[maxbibnames=99]{biblatex}
\usepackage{fancybox,graphicx}
\usepackage[title]{appendix}
\usepackage{xcolor}
\usepackage{url}
\usepackage[linkbordercolor={1 0 0}]{hyperref}
\usepackage{mathtools}
\usepackage{amssymb}
\usepackage{algorithm}
\usepackage[noend]{algpseudocode}
\usepackage{authblk}
\DeclarePairedDelimiter{\ceil}{\lceil}{\rceil}

\hide{
\usepackage{amscd}
\usepackage{amsfonts}
\usepackage{amsmath}
\usepackage{amssymb}
\usepackage{amsthm}
\usepackage{cases}		
\usepackage{cutwin}
\usepackage{enumerate}
\usepackage{epstopdf}
\usepackage{graphicx}
\usepackage{ifthen}
\usepackage{lipsum}
\usepackage{mathrsfs}	
\usepackage{multimedia}
\usepackage{wrapfig}
}
\usepackage{amsthm}
\newtheorem{theorem}{Theorem}[section]

\newtheorem{lemma}[theorem]{Lemma}
\newtheorem{proposition}[theorem]{Proposition}

\newlength\tindent
\title{ New First-Order Algorithms for Stochastic Variational Inequalities }

\author{
Kevin Huang\thanks{Department of Industrial and System Engineering, University of Minnesota, huan1741@umn.edu}
\hspace{1cm}
Shuzhong Zhang\thanks{Department of Industrial and System Engineering, University of Minnesota, zhangs@umn.edu}
}

\date{\today}

\begin{document}

\maketitle

\begin{abstract}
    In this paper, we propose two new solution schemes to solve the stochastic strongly monotone variational inequality problems: the stochastic extra-point solution scheme and the stochastic extra-momentum solution scheme. The first one is a general scheme based on updating the iterative sequence and an auxiliary extra-point sequence. In the case of deterministic VI model, this approach includes several state-of-the-art first-order methods as its special cases. The second scheme combines two momentum-based directions: the so-called {\it heavy-ball}\/ direction and the {\it optimism}\/ direction, where only one projection per iteration is required in its updating process. We show that, if the variance of the stochastic oracle is appropriately controlled, then both schemes can be made to achieve optimal iteration complexity of $\mathcal{O}\left(\kappa\ln\left(\frac{1}{\epsilon}\right)\right)$ to reach an $\epsilon$-solution for a strongly monotone VI problem with condition number $\kappa$. We show that these methods can be readily incorporated in a zeroth-order approach to solve stochastic minimax saddle-point problems, where only noisy and biased samples of the objective can be obtained, with a total sample complexity of $\mathcal{O}\left(\frac{\kappa^2}{\epsilon}\ln\left(\frac{1}{\epsilon}\right)\right)$.

    \vspace{3mm}

    \noindent\textbf{Keywords:} variational inequality, minimax saddle-point, stochastic first-order method, zeroth-order method. %stochastic zeroth-order methods.
\end{abstract}

\section{Introduction}
\label{sec:intro}

Given a constraint set $\mathcal{Z}\subseteq \mathbb{R}^n$ and a mapping $F:\mathcal{Z}\rightarrow \mathbb{R}^n$, the classical variational inequality (VI) problem is to find $z^*\in\mathcal{Z}$ such that
\begin{eqnarray}
    &&F(z^*)^{\top}(z-z^*)\geq 0,\quad\forall z\in\mathcal{Z}.\label{vi-prob}
\end{eqnarray}
%The study of VI problem \eqref{vi-prob} dates back to 1960's, first in the form of complementarity problem (CP) which models various equilibrium settings such as the economic supply-demand equilibria, traffic equilibria, and generally the Nash equilibria.
For an introduction to VI and its applications, we refer the readers to
Facchinei and Pang~\cite{facchinei2007finite} and the references therein.

In this paper, we %take one step further and
consider a {\it stochastic}\/ version of problem \eqref{vi-prob}, where the exact evaluation of the mapping $F(\cdot)$ is inaccessible. Instead, only a {\it stochastic oracle}\/ is available.  The stochasticity in question may stem from, e.g., the non-deterministic nature of mixed strategies of the players in a game-setting, or simply because of the difficulty in evaluating the mapping itself. The latter has become more pronounced in the literature, due to its recent-found application as a training/learning subproblem in machine learning and/or statistical learning.
%In fact, VI includes minimax saddle-point problems as a special case, which is capable of modeling, for instance, the generative adversarial network (GAN) \cite{goodfellow2014generative, gidel2018variational}.
%In order to generate the update direction in our proposed methods, we require the access to the
The so-called {\it stochastic oracle}\/ is a noisy estimation of the mapping $F(\cdot)$, and an iterative scheme that incorporates such oracle is known as {\it stochastic approximation}\/ (SA).
%We shall reserve the detailed assumptions for such oracle as well as monotonicity until Section \ref{sec:sto-VI-pre}.
%As a matter of fact, the idea of using
As far as we know, the first proposal to use such approach for stochastic optimization can be traced back to the seminal work of %was first proposed and studied by
Robbins and Monro \cite{robbins1951stochastic}.
%However, in contrast to the early developments of stochastic optimization, the research on VI has been focusing on the deterministic settings and it is not until
In 2008, Jiang and Xu \cite{jiang2008stochastic} followed the SA approach %introduced this idea
to solve VI models. Since then, efforts have been made to extend existing deterministic methods to the stochastic VI models; see e.g.~\cite{juditsky2011solving, yousefian2014optimal, kannan2019optimal, kotsalis2020simple, jalilzadeh2019proximal, iusem2017extragradient, iusem2019variance}.

Let us start our discussion by introducing the assumptions made throughout the paper. We consider VI model~\eqref{vi-prob} where $\mathcal{Z}$ is a closed convex set. Moreover, the following two conditions are assumed:
\begin{equation}
    \label{pre:strong-mono}
    \left(F(z)-F(z')\right)^{\top}(z-z')\geq \mu\|z-z'\|^2, \quad\forall z,z'\in\mathcal{Z},
\end{equation}
for some $\mu>0$, and
\begin{equation}
    \label{pre:lip}
    \|F(z)-F(z')\|\leq L\|z-z'\|,\quad\forall z,z'\in \mathcal{Z},
\end{equation}
for some $L \ge \mu>0$. Condition~\eqref{pre:strong-mono} is known as the {\it strong monotonicity}\/ of $F$, while Condition~\eqref{pre:lip} is known as the {\it Lipschitz continuity}\/ of $F$. If Condition \ref{pre:strong-mono} is met with $\mu=0$ then $F$ is known as {\it monotone}. VI problems satisfying \eqref{pre:strong-mono} with positive $\mu$
can be easily shown to have a unique solution $z^*$. Let us denote
    $\kappa := \frac{L}{\mu}\ge 1$.
Parameter $\kappa$ is usually known as the condition number of the VI model~\eqref{vi-prob}.
We also assume
\begin{eqnarray}
    &&\max\limits_{z,z'\in\mathcal{Z}}\|z-z'\|\le D,\label{feasible-set-bound}
\end{eqnarray}
namely the constraint set is bounded. Remark that this assumption can actually be removed without affecting the results, but then the analysis becomes lengthier and tedious without conceivable conceptual benefit, and so we shall not pursue that generality in this paper.
%an increased level of technicalities.
%However, we shall remark that this assumption is purely for the purpose of simplifying the analysis, but not required for our methods. It can be eliminated through proper (but more involved) analysis.

The {\it stochastic oracle}\/ of the mapping, denote by $\hat{F}(z,\xi)$, takes a random sample $\xi\in\Xi$ from some sample space $\Xi$. The oracle is required to satisfy:
\begin{eqnarray}
    \mathbb{E}\left[\|\hat{F}(z,\xi)-F(z)\|\right]&\le& \delta\label{oracle-bias}\\
    \mathbb{E}\left[\|\hat{F}(z,\xi)-F(z)\|^2\right]&\le& \sigma^2\label{oracle-var}
\end{eqnarray}
for all $z\in\mathcal{Z}$, where $\delta,\sigma\ge 0$ are some constants. In other words, we assume both the bias and the deviation are uniformly upper-bounded.

In this paper, we propose two stochastic first-order schemes: the {\it stochastic extra-point scheme}\/ and the {\it stochastic extra-momentum scheme}. The first scheme maintains two sequences of iterates %throughout the iterations
%and features combination of
featuring several well-known first-order search directions %extending several 
such as: the extra-gradient \cite{korpelevich1976extragradient, tseng1995linear}, the heavy-ball \cite{polyak1964some}, Nesterov's extrapolation \cite{nesterov1983method, nesterov2003introductory}, and the optimism direction~\cite{popov1980modification, mokhtari2019unified}. The second scheme, on the other hand, specifically combines the {\it heavy-ball}\/ momentum and the {\it optimism}\/ momentum in its updating formula, and maintains only one sequence throughout the iterations, therefore requiring only {\it one}\/ projection per iteration.
%A different treatment in the analysis is established for the second scheme in order to show optimal iteration complexity.
These two approaches require different types of analysis.
Both schemes render a wider range of search directions than the existing first-order methods, and the parameters associated with each search direction could and {\it should}\/ be tuned differently from problem-class to problem-class in order to secure good practical performances. The deterministic counterpart of these methods can be found in our previous work \cite{huang2021unifying}. In the stochastic context, we show that as long as the variance can be reduced throughout the iterations, they yield the {\it optimal}\/ iteration complexity 
$\mathcal{O}(\kappa\ln(1/\epsilon))$ (cf.~\cite{junyu2019}) to reach $\epsilon$-solution: $\|z^k-z^*\|^2\le\epsilon$, with an additional biased term depending on $\delta$. In a later section, we demonstrate an application to the stochastic black-box minimax saddle-point problem %(aka the minimax problem)
where only noisy function values $f(x,y)$ are accessible. This application is particularly relevant, given its applications in machine learning, where the training data set may be very large and evaluating exact gradient/function value is usually impractical.
%could worsen the performance.
Through a smoothing technique, we propose a {\it stochastic zeroth-order gradient}\/ as our update directions in either the stochastic extra-point scheme or the stochastic extra-momentum scheme. We show that both approaches yield an iteration complexity of $\mathcal{O}(\kappa\ln(1/\epsilon))$ and a sample-complexity of $\mathcal{O}\left(\frac{\kappa^2}{\epsilon}\ln\left(\frac{1}{\epsilon}\right)\right)$.
%. Iteration complexity of that approach is studied, together with a sample-complexity analysis. % and reduce the variance through mini-batch sampling.

The rest of the paper is organized as follows. In Section \ref{sec:literature}, we survey the relevant literature with a focus on stochastic VI. In Section \ref{sec:sto-VI}, we present the main results in this paper, i.e.\ the convergence results of the two proposed methods, while the technical proofs are relegated to the appendices. In Section \ref{sec:sto-0-minmax}, we introduce a stochastic black-box saddle-point problem and present the sample complexity results of our methods. We present some promising preliminary numerical results in Section \ref{sec:numerical} and conclude the paper in Section \ref{sec:conclusion}.

\section{Literature Review}
\label{sec:literature}

The first-order algorithms for deterministic VI \eqref{vi-prob} serve as a basis for the developments of their stochastic counterparts. These algorithms include the projection method~\cite{facchinei2007finite}, the proximal method \cite{martinet1970breve, rockafellar1976monotone, tseng1995linear}, the extra-gradient method \cite{korpelevich1976extragradient, tseng1995linear}, the optimistic gradient descent ascent (OGDA) method \cite{popov1980modification, mokhtari2019unified, mokhtari2020convergence}, the mirror-prox method \cite{nemirovski2004prox}, the extrapolation method \cite{nesterov2006solving, nesterov2007dual}, and the extra-point method~\cite{huang2021unifying}.
%For lower bound analysis, see \cite{junyu2019}.

In this section, we shall focus on the developments of algorithms for stochastic VI, starting with a paper of Jiang and Xu \cite{jiang2008stochastic}, where the authors propose a stochastic projection method for solving strongly monotone and Lipschitz continuous VI problems and present an almost-sure convergence result. Koshal {\it et al}.\ \cite{koshal2012regularized} propose iterative Tikhonov regularization method and iterative proximal point method and show almost-sure convergence with monotone and Lipschitz continuous VI problems. Both methods solve a strongly monotone VI subproblem at each iteration. Yousefian {\it et al}.\ \cite{yousefian2013regularized} further introduce local smoothing technique to the above-mentioned regularized methods to account for non-Lipschitz mappings and show almost-sure convergence. A survey on these methods, as well as applications and the theory behind stochastic VI can be found in Shanbhag \cite{shanbhag2013stochastic}.

Juditsky {\it et al}.\ \cite{juditsky2011solving} are among the first to show an iteration complexity bound for stochastic VI algorithms. They extend the mirror-prox method \cite{nemirovski2004prox} to stochastic settings and prove an optimal iteration complexity bound for monotone VI: $\mathcal{O}(\frac{1}{\epsilon^2})$ , or $\mathcal{O}\left(\frac{1}{\epsilon}\right)$ when the variance can be controlled small enough. Yousefian {\it et al}.\ \cite{yousefian2014optimal} further extend the stochastic mirror-prox method with a more general step size choice and show an $\mathcal{O}\left(\frac{1}{\epsilon^2}\right)$ iteration complexity, where they also show an $\mathcal{O}(\frac{1}{\epsilon})$ complexity for the stochastic extra-gradient method for solving strongly monotone VI problems. Yousefian {\it et al.} \cite{yousefian2017smoothing} use randomized smoothing technique for non-Lipschitz mapping and show an $\mathcal{O}\left(\frac{1}{\epsilon^6}\right)$ iteration complexity. Chen {\it et al}.\ \cite{chen2017accelerated} consider a specific class of VI model: a mapping that consists of a Lipschitz continuous and monotone operator, a Lipschitz continuous gradient mapping of a convex function, and a subgradient mapping of a simple convex function. They propose a method that combines Nesterov's acceleration \cite{nesterov2003introductory} with the stochastic mirror-prox method to exploit this special structure, resulting in an optimal iteration complexity for such class of problem: $\mathcal{O}\left(\frac{1}{\epsilon^2}\right)$, or $\mathcal{O}\left(\frac{1}{\epsilon}\right)$ when the variance can be controlled small enough, or $\mathcal{O}\left(\sqrt{\frac{1}{\epsilon}}\right)$ when the operator consists only of gradient/subgradient mappings from some convex function. Kannan and Shanbhag \cite{kannan2019optimal} analyze a general variant of extra-gradient method (which uses general distance-generating functions) and show that under a slightly weaker assumptions than the strongly monotonicity, the optimal $\mathcal{O}\left(\frac{1}{\epsilon}\right)$ iteration bound still hold. Kotsalis {\it et al.} \cite{kotsalis2020simple} extend the OGDA method to strongly monotone stochastic VI with iteration complexity $\mathcal{O}\left(\frac{1}{\epsilon}\right)$, or $\mathcal{O}\left(\kappa\ln\left(\frac{1}{\epsilon}\right)\right)$ when the variance can be controlled small enough.

We shall note that the detailed implementation of variance-reduction is in general not considered in the above-mentioned methods (although some do present additional complexity term when the variance is small, such as in \cite{juditsky2011solving, chen2017accelerated}). Therefore, the optimal iteration complexity bound is $\mathcal{O}\left(\frac{1}{\epsilon^2}\right)$ for monotone VI and $\mathcal{O}\left(\frac{1}{\epsilon}\right)$ for strongly monotone VI, as compared to $\mathcal{O}\left(\frac{1}{\epsilon}\right)$ and $\mathcal{O}(\kappa\ln(1/\epsilon))$ for their deterministic counterpart. By increasing the sample size (aka mini-batch) in each iteration, the variance can be reduced as the algorithm progresses, therefore attaining the same optimal iteration complexity bound as the deterministic problems. %A bound with multiple terms involved indicates a combination of worst cast complexity bounds when variance is controlled and when it is not.

There have been developments for variance-reduction-based methods in recent years. Jalilzadeh and Shanbhag \cite{jalilzadeh2019proximal} extend the method \cite{nesterov2006solving} for deterministic strongly monotone VI to stochastic VI and show that with variance reduction the optimal iteration complexity $\mathcal{O}(\kappa\ln(1/\epsilon))$ can be achieved, together with a total sample complexity of $\mathcal{O}\left(\frac{1}{\epsilon^{\beta}}\right)$ for some constant $\beta>1$. With this method as a subroutine, they also propose a variance-reduced proximal point method with iteration complexity $\mathcal{O}\left(\frac{1}{\epsilon}\ln\left(\frac{1}{\epsilon}\right)\right)$ and sample complexity $\mathcal{O}(\frac{1}{\epsilon^{1+2\alpha\beta}})$ for some constants $\alpha,\beta>1$. Iusem {\it et al.}\/ \cite{iusem2017extragradient} propose a variance-reduced extra-gradient-based method for monotone VI and show $\mathcal{O}\left(\frac{1}{\epsilon}\right)$ iteration complexity and $\mathcal{O}\left(\frac{1}{\epsilon}\right)$ sample complexity. They further extend the method \cite{iusem2019variance} by incorporating line-search for unknown Lipschitz constant, while preserving similar bounds. Palaniappan and Bach \cite{palaniappan2016stochastic} propose variance-reduced stochastic forward-backward methods based on (accelerated) stochastic gradient descent methods in optimization and show $\mathcal{O}(\kappa\ln(1/\epsilon))$ iteration complexity. For another line of work, which includes the concept of differential privacy in the stochastic VI, we refer the readers to a recent paper \cite{boob2021optimal} and the references therein.
The stochastic oracle maybe {\it man-made}. For instance, the technique of randomized smoothing has been applied in the so-called zeroth-order methods (i.e.\ derivative-free methods), refer to \cite{nesterov2017random, shalev2011online} or the survey \cite{larson2019derivative} in the context of optimization and \cite{wang2020zeroth,xu2020gradient,liu2019min,roy2019online,menickelly2020derivative} in the context of minimax saddle-point problems.

\section{The Stochastic First-Order Methods for Strongly Monotone VI}
\label{sec:sto-VI}

Let us start this section by introducing the notations to facilitate our analysis.  %introducing a simplification of notation to facilitate the presentation of the analysis.
We shall denote the stochastic oracle as $\hat{F}(\cdot)$, suppressing the notation $\xi$ whenever it is clear from the context. For example, $\hat{F}(z^k)$ is associated with the random sample $\xi^k\in\Xi$. In addition, we denote $P_{\mathcal{Z}}(\cdot)$ as the {\it projection operator} onto the feasible set $\mathcal{Z}$.

\subsection{The stochastic extra-point scheme}
\label{sec:sto-exp}

We first present the iterative updating rule for the stochastic extra-point scheme:
\begin{eqnarray}
\left\{
\begin{array}{lcl}
     z^{k+0.5} &:=& P_{\mathcal{Z}}\left(z^k+\beta(z^k-z^{k-1})-\eta \hat{F}(z^k)\right), \\
     z^{k+1}   &:=& P_{\mathcal{Z}}\left(z^k-\alpha \hat{F}(z^{k+0.5})+\gamma(z^k-z^{k-1})-\tau(\hat{F}(z^{k})-\hat{F}(z^{k-1}))\right), \\
\end{array}
\right.\label{sto-exp-update}
\end{eqnarray}
for $k=0,1,...$, where the sequence $\{ z^k \mid k=0,1,...\}$ is the sequence of iterates, and $\{ z^{k+0.5} \mid k=0,1,...\}$ is the sequence of {\it extra points}, which helps to produce the sequence of iterates.

In the case of {\it deterministic}\/ strongly monotone VI, we introduced in our previous work \cite{huang2021unifying} a {\it unifying}\/ extra-point updating scheme, %\eqref{sto-exp-update}
which includes specific first-order search directions such as the extra-gradient, the heavy-ball method, the optimistic method, and Nesterov's extrapolation; these are incorporated with the parameters $\alpha,\beta,\gamma,\eta,\tau\ge0$. As any specific configuration of these parameters should be tailored to the problem structure at hand, our goal is to provide conditions of the parameters under which an {\it optimal}\/ iteration complexity can be guaranteed. %A simple example that satisfies these conditions will also be provided.
This line of analysis will now be extended to solve stochastic VI as given in \eqref{sto-exp-update}.
We shall first establish the relational inequalities between subsequent iterates in terms of the {\it expected}\/ distance to the unique solution $z^*$, denote by $d_k = \mathbb{E}\left[\|z^k-z^*\|^2\right]$.

\begin{lemma}
\label{lem:sto-exp-relation}
For the sequences $\{z^k\mid k=0,1,...\}$ and $\{z^{k+0.5}\mid k=0,1,...\}$ generated from the stochastic extra-point scheme \eqref{sto-exp-update}, the following inequality holds:
\begin{eqnarray}
&&(1-4|\gamma-\beta|-\tau L)d_{k+1}\nonumber\\
&\le& \left(1+4\gamma+6|\gamma-\beta|+4\tau L-\alpha\mu\right)d_k %\nonumber\\
+\left(2|\gamma-\beta|+2\gamma+4\tau L\right)d_{k-1}\nonumber\\
&&+\left(2\alpha^2L^2+2|\gamma-\beta|+2\gamma+2\alpha\mu-1\right)\mathbb{E}\left[\|z^{k+0.5}-z^k\|^2\right]\nonumber\\
&&+8\left(\alpha^2+\frac{\tau}{L}\right)\sigma^2+2\alpha\delta D %\nonumber\\
+2(\eta-\alpha)\mathbb{E}\left[(z^{k+1}-z^{k+0.5})^\top\hat{F}(z^k)\right].\label{sto-exp-relation}
\end{eqnarray}
\end{lemma}
\begin{proof}
See Appendix \ref{proof:sto-exp-relation}.
\end{proof}

Lemma \ref{lem:sto-exp-relation} forms a basis to the desired linear convergence, and it is possible to identify the conditions for the parameters $\alpha,\beta,\gamma,\eta,\tau$ in order to achieve linear convergence. Consider parameters satisfying %the following two constraints:
\begin{eqnarray}
\left\{
\begin{array}{ll}
     \eta=\alpha,\\
     2\alpha^2L^2+2|\gamma-\beta|+2\gamma+2\alpha\mu-1 \le 0,
\end{array}
\right.\label{sto-exp-cond-1}
\end{eqnarray}
and denote
\begin{eqnarray}
    \left\{
    \begin{array}{ll}
         t_1 = \alpha\mu-4\gamma-6|\gamma-\beta|-4\tau L, \\
         t_2 = 2|\gamma-\beta|+2\gamma+4\tau L,\\
         t_3 = 4|\gamma-\beta|+\tau L.
    \end{array}
    \right.\label{sto-exp-t}
\end{eqnarray}
Then we obtain from \eqref{sto-exp-relation} that
\begin{eqnarray}
(1-t_3)d_{k+1}&\le& \left(1-t_1\right)d_k+t_2d_{k-1}+8\left(\alpha^2+\frac{\tau}{L}\right)\sigma^2+2\alpha\delta D.\label{sto-exp-relation-2}
\end{eqnarray}
%The above inequality \eqref{sto-exp-relation-2} is one step away from our final result.
With additional constraints on $t_1,t_2,t_3$, the {\it variance-reduced} convergence result is summarized in the next theorem.

\begin{theorem}
\label{th:sto-exp-result}
For non-negative parameters $\alpha,\beta,\gamma,\eta,\tau$ satisfying \eqref{sto-exp-cond-1} and \eqref{sto-exp-t}, suppose that
    \begin{eqnarray}
    \left\{
    \begin{array}{ll}
         0\le t_3 < t_1<1, \\
         t_2 < t_1-t_3.
    \end{array}
    \right.\label{sto-exp-cond-t}
    \end{eqnarray}
Let $q=\frac{2(1-t_3)}{t_1-t_2-t_3}>1$. For a fixed precision $\epsilon>0$, denote $K=\mathcal{O}\left(q\cdot\ln\left(\frac{1}{\epsilon}\right)\right)$.
%Suppose furthermore that the variance is controlled so that
%    \(
%    \sigma^2= \mathcal{O} \left(\left(1-\frac{1}{q}\right)^K\right) = \mathcal{O} (\epsilon).
%    \)
%For a fixed iteration count $K$ and the sequence $\{z^k\}$, $k=0,1,..,K$ generated from the stochastic extra-point method \eqref{sto-exp-update}, the expected distance to the solution of iterate $z^K$ is bounded by:
Then, we have
\begin{eqnarray}
d_K=\mathbb{E}\left[\|z^K-z^*\|^2\right]\le \mathcal{O}(\epsilon) + \mathcal{O}(\sigma^2) + \mathcal{O}(\delta D). \label{sto-exp-final-bound}
\end{eqnarray}
%with $K=\mathcal{O}\left(q\cdot\ln\left(\frac{1}{\epsilon}\right)\right)$, $q=\frac{2(1-t_3)}{t_1-t_2-t_3}>1$, given the following two conditions are satisfied:
%\begin{enumerate}
%    \item The non-negative parameters $\alpha,\beta,\gamma,\eta,\tau$ satisfy \eqref{sto-exp-cond-1}, \eqref{sto-exp-t}, and additionally
%    \begin{eqnarray}
%    \left\{
%    \begin{array}{ll}
%         0\le t_3 < t_1<1, \\
%         t_2 < t_1-t_3.
%    \end{array}
%    \right.\label{sto-exp-cond-t}
%    \end{eqnarray}
%    \item The variance is controlled such that
%    \[
%    \sigma^2={\color{red}o}\left(\left(1-\frac{1}{q}\right)^K\right) {\color{red} = O(\epsilon)}.
%    \]
%    \end{enumerate}
%{\color{red} (Why small o? It should be big O. Right? Is $\sigma^2=O(\epsilon)$?}
\end{theorem}

\begin{proof}
See Appendix \ref{proof:sto-exp-result}.
\end{proof}

Regarding Theorem \ref{th:sto-exp-result}, a few remarks are in order.
%Firstly, there is a bias term of the order $\mathcal{O}(\delta D)$ present in \eqref{sto-exp-final-bound} because we only assume the bias of the stochastic oracle $\hat{F}(\cdot)$ is bounded by $\delta\ge0$. For an unbiased oracle, namely $\delta=0$, the bias term will no longer present in our result.
First, as we remarked earlier, the boundedness condition \eqref{feasible-set-bound} can be removed. However, the analysis will become much longer and cumbersome; we keep it here for simplicity.
%is introduced, such bound is merely for the convenience in the derivation and is by no means relevant in our method. However, we shall omit the more involved proof for space concern.
%Thirdly, we {\it assumed} the variance can be controlled to have a uniform upper bound in the second condition of Theorem \ref{th:sto-exp-result} without specifying {\it how} to control it.
Second, a common way to achieve variance reduction is through increasing the mini-batch sample sizes. In fact, we may fix the sample size at the beginning with order $\left(1-\frac{1}{q}\right)^{-K}$, or it increases linearly at a rate $\left(1+\frac{1}{q}\right)$ as $k$ increases. We shall discuss more on this strategy in Section \ref{sec:sto-0-minmax}. % for an application in stochastic black-box saddle point problems.
Finally, we note that without variance reduction, it is possible to adopt {\it diminishing}\/ step sizes $(\alpha_k,\beta_k,\gamma_k,\eta_k,\tau_k)$ instead of fixing step sizes as we have assumed so far. The {\it optimal}\/ uniform sublinear convergence rate $\frac{1}{k}$ can be established through a separate analysis continued from Lemma \ref{lem:sto-exp-relation}.
%The interested readers could find the detailed analysis in
The details can be found in Appendix \ref{appendix:sublinear-conv}.

Next proposition concludes this subsection with a specific choice of the parameters.
%This example illustrates how the proposed method yields an optimal iteration complexity $\mathcal{O}(\kappa\ln(1/\epsilon))$, with a presumably controllable variance.

\begin{proposition}
\label{prop:sto-exp-example}
If one chooses the following parameters
\begin{eqnarray}
(\alpha,\beta,\gamma,\eta,\tau)=\left(\frac{1}{4L},\frac{1}{64\kappa},\frac{1}{64\kappa},\frac{1}{4L},\frac{1}{64L\kappa}\right)\nonumber
\end{eqnarray}
in \eqref{sto-exp-update} (thus $(t_1,t_2,t_3)=\left(\frac{1}{8\kappa},\frac{3}{32\kappa},\frac{1}{64\kappa}\right)$) then it holds that
\[
d_k\le \left(1-\frac{1}{256\kappa}\right)^{k}\cdot\frac{283}{256}\|z^0-z^*\|^2+\left(\frac{40\sigma^2}{63L^2}+\frac{32\delta D}{63L}\right)\cdot256\kappa.
\]
\end{proposition}
\begin{proof}
See Appendix \ref{proof:sto-exp-example}.
\end{proof}

\subsection{The stochastic extra-momentum scheme}
\label{sec:sto-exm}

In this subsection, we present an alternative stochastic first-order method that achieves the optimal iteration complexity as well, the {\it stochastic extra-momentum scheme}:
\begin{eqnarray}
z^{k+1} := P_{\mathcal{Z}}\left(z^k-\alpha \hat{F}(z^k)+\gamma(z^k-z^{k-1})-\tau\left(\hat{F}(z^k)-\hat{F}(z^{k-1})\right)\right),\label{sto-exm-update}
\end{eqnarray}
for $k=0,1,...$.

Compared with the stochastic extra-point scheme~\eqref{sto-exp-update}, the above update \eqref{sto-exm-update} manipulates only %exhibits limitations on the generalization as it only adopts
the {\it momentum}\/ terms alongside the stochastic gradient direction (the notion ``gradient'' here refers to the mapping $F(\cdot)$ in the VI model), namely the heavy-ball direction $z^k-z^{k-1}$ and the optimism direction $\hat{F}(z^k)-\hat{F}(z^{k-1})$. %in its search directions.
%However, the advantage is also obvious:
Since it maintains a single sequence $\{z^k\}$ throughout the iterations, this scheme requires {\it one}\/ projection per iteration, as compared to two projections in the case of the stochastic extra-point scheme. %In fact, as far as we know, most state-of-the-art stochastic first-order methods that achieve the {\it optimal}\/ iteration complexity.
We shall remark that the method proposed in Kotsalis {\it et al.}~\cite{kotsalis2020simple} only considers the optimism term. Therefore, the stochastic extra-momentum scheme introduced above may be viewed as a generalization.

As in the previous subsection, we shall first establish a relational inequality between the iterates.
%, similar to what we did in the previous section.  However,
As we can see from the lemma below, the structure of this relational inequality is in fact quite different from the previous case. The detailed proof can be found in the appendix.

\begin{lemma}
\label{lem:sto-exm-relation}
For the sequence $\{z^k\mid k=0,1,...\}$ generated from the stochastic extra-momentum scheme \eqref{sto-exm-update}, the following inequality holds:
\begin{eqnarray}
&&\mathbb{E}\left[\left(\frac{1}{2}+\frac{\alpha\mu}{2}-\frac{\gamma}{2}\right)\|z^{k+1}-z^*\|^2+\alpha(z^{k+1}-z^*)^\top\left(\hat{F}(z^k)-\hat{F}(z^{k+1})\right)+\frac{1}{4}\|z^{k+1}-z^k\|^2\right]\nonumber\\
&\le& \mathbb{E}\left[\frac{1}{2}\|z^k-z^*\|^2+\tau(z^k-z^*)^\top\left(\hat{F}(z^{k-1})-\hat{F}(z^{k})\right)+\left(2\tau^2L^2+\frac{\gamma}{2}\right)\|z^k-z^{k-1}\|^2\right]\nonumber\\
&&+8\tau^2\sigma^2+\frac{\alpha\delta^2}{2\mu}.\label{sto-exm-relation}
\end{eqnarray}
\end{lemma}
\begin{proof}
See Appendix \ref{proof:sto-exm-relation}.
\end{proof}

%In view of the structure in \eqref{sto-exm-relation}, we
Observe that each of the terms on the LHS of \eqref{sto-exm-relation} differs in the iteration index from the RHS exactly by one.
%, in terms of the iteration index $k$, from each of the terms on the RHS.
This property enables us to design a possible {\it potential function}\/ that measures the convergence of the iterative process. %, which we shall specify shortly after.
We shall specify additional conditions on the non-negative parameters $\alpha,\gamma,\tau$ in order to further simplify \eqref{sto-exm-relation}:
\begin{eqnarray}
&& 1+\alpha\mu-\gamma\ge 1+\frac{\theta}{\kappa},\quad \frac{\alpha}{\tau}= 1+\frac{\theta}{\kappa},\quad \frac{1}{8\tau^2L^2+2\gamma}\ge 1+\frac{\theta}{\kappa},\label{sto-exm-condition-par}
\end{eqnarray}
where $\theta\in(0,1]$ is some constant independent of $\kappa$. Note that the LHS of each inequality in \eqref{sto-exm-condition-par} is the ratio between the coefficients on the LHS and RHS of \eqref{sto-exm-relation} for each corresponding term. Therefore, the relation~\eqref{sto-exm-relation} can be rearranged as:
\begin{eqnarray}
&& \left(1+\frac{\theta}{\kappa}\right)\mathbb{E}\left[\frac{1}{2}\|z^{k+1}-z^*\|^2+\tau(z^{k+1}-z^*)^\top\left(\hat{F}(z^{k})-\hat{F}(z^{k+1})\right)+\left(2\tau^2L^2+\frac{\gamma}{2}\right)\|z^{k+1}-z^{k}\|^2\right]\nonumber\\
&\le&\mathbb{E}\left[\left(\frac{1}{2}+\frac{\alpha\mu}{2}-\frac{\gamma}{2}\right)\|z^{k+1}-z^*\|^2+\alpha(z^{k+1}-z^*)^\top\left(\hat{F}(z^k)-\hat{F}(z^{k+1})\right)+\frac{1}{4}\|z^{k+1}-z^k\|^2\right]\nonumber\\
&\le& \mathbb{E}\left[\frac{1}{2}\|z^k-z^*\|^2+\tau(z^k-z^*)^\top\left(\hat{F}(z^{k-1})-\hat{F}(z^{k})\right)+\left(2\tau^2L^2+\frac{\gamma}{2}\right)\|z^k-z^{k-1}\|^2\right]\nonumber\\
&&+8\tau^2\sigma^2+\frac{\alpha\delta^2}{2\mu}.\label{sto-exm-relation-2}
\end{eqnarray}

Now, by defining the potential function $V_k$ as
\[
V_k = \mathbb{E}\left[\frac{1}{2}\|z^k-z^*\|^2+\tau(z^k-z^*)^\top\left(\hat{F}(z^{k-1})-\hat{F}(z^{k})\right)+\left(2\tau^2L^2+\frac{\gamma}{2}\right)\|z^k-z^{k-1}\|^2\right],
\]
inequality \eqref{sto-exm-relation-2} can be rewritten as
\begin{eqnarray}
\left(1+\frac{\theta}{\kappa}\right)V_{k+1}\le V_k +8\tau^2\sigma^2+\frac{\alpha\delta^2}{2\mu}.\label{sto-exm-V-relation}
\end{eqnarray}
This leads to our final results, as summarized in the next theorem:

\begin{theorem}
\label{th:sto-exm-result}
Suppose that the non-negative parameters $\alpha,\gamma,\tau$ satisfy \eqref{sto-exm-condition-par} for some constant $\theta\in(0,1]$.
For %any fixed iteration count $K$ and
the sequence $\{z^k \mid k=0,1,...\}$
generated from the stochastic extra-momentum scheme \eqref{sto-exm-update}, the expected distance to the solution of iterate $z^k$ is bounded by:
\begin{eqnarray}
%d_K=\mathbb{E}\left[\|z^K-z^*\|^2\right]\le \epsilon+\mathcal{O}(\delta ^2).
\mathbb{E}\left[\|z^k-z^*\|^2\right] \le
2\left(1+\frac{\theta}{\kappa}\right)^{-k}\|z^0-z^*\|^2+\left(\frac{\kappa}{\theta}+1\right)\cdot 32\tau^2\sigma^2+\frac{2\kappa\alpha\delta^2}{\theta\mu}.
\label{sto-exm-final-bound}
\end{eqnarray}
%with $K=\mathcal{O}\left(\kappa\ln\left(\frac{1}{\epsilon}\right)\right)$, given the following two conditions are satisfied:
%\begin{enumerate}
%    \item The non-negative parameters $\alpha,\gamma,\tau$ satisfy \eqref{sto-exm-condition-par} for some constant $\theta\in(0,1]$.
%    \item The variance is controlled such that
%    \[
%    \sigma^2=o\left(\left(1-\frac{1}{\kappa}\right)^K\right).
%    \]
%\end{enumerate}
\end{theorem}
\begin{proof}
See Appendix \ref{proof:sto-exm-result}.
\end{proof}

%We remark that the bound in \eqref{sto-exm-final-bound} is biased by $\mathcal{O}(\delta^2)$ due to our assumption \eqref{oracle-bias}, and that the uniform bound on the variance $\sigma^2$ in the second condition can be achieved by increasing the sample sizes linearly throughout the iterations, which we will discuss in the next section.
%
%We present a simple example of a parameter choice satisfying \eqref{sto-exm-condition-par} in the next proposition to conclude this section.

A simple choice of parameters leads to:
\begin{proposition}
\label{prop:sto-exm-example}
If we choose the parameters as
\[
(\alpha,\tau,\gamma)=\left(\frac{1}{4L},\frac{\alpha}{1+\frac{\theta}{\kappa}},\frac{1}{8(\kappa+\theta)}\right),\quad \theta=\frac{1}{8},
\]
then scheme~\eqref{sto-exm-update} assures that
\begin{eqnarray}
&& d_k \le 2\left(1-\frac{1}{8\kappa+1}\right)^k\|z^0-z^*\|^2+\frac{128\sigma^2}{\mu(8L+\mu)}+\frac{4\delta^2}{\mu^2}.
\end{eqnarray}
\end{proposition}
\begin{proof}
It follows by substituting the parameter choice into \eqref{sto-exm-final-bound}.
\end{proof}

This shows that if we run the stochastic extra-momentum scheme~\eqref{sto-exm-update} with the above parameter choice, then in $K=\mathcal{O}(\kappa \ln \frac{1}{\epsilon})$ iterations we will reach a solution satisfying
\[
\mathbb{E}\left[\|z^K-z^*\|^2\right] \le \mathcal{O}(\epsilon) + \mathcal{O} \left( \frac{\sigma^2}{\mu L}\right) + \mathcal{O} \left( \frac{\delta^2}{\mu^2} \right).
\]

%{\color{blue} (Shall we present an example of total sample complexity, to be compared to the results in the literature, say by Uday's group?)}

\section{A Stochastic Zeroth-Order Approach to Saddle-Point Problems}
\label{sec:sto-0-minmax}
%\subsection{Preparations} \label{sec:prep}
In this section, we shall apply the proposed stochastic extra-point/extra-momentum scheme to solve the following saddle-point problem without needing to compute the gradients of $f$:
\begin{eqnarray}
\min\limits_{x\in\mathcal{X}}\max\limits_{y\in\mathcal{Y}}f(x,y)\label{spp},
\end{eqnarray}
where $\mathcal{X}\subseteq \mathbb{R}^n$, $\mathcal{Y}\subseteq\mathbb{R}^m$ are convex sets, $f(x,y)$ is strongly convex (with fixed $y$) and strongly concave (with fixed $x$) with modulus $\mu$, and the partial gradients $\nabla_xf(x,y)$/$\nabla_yf(x,y)$ are Lipschitz continuous with constant $L_x/L_y$ for fixed $y/x$, and with constant $L_{xy}$ with fixed $x/y$. We let $L=2 \cdot\max(L_x,L_y,L_{xy})$. We further assume that the function $f(x,y)$ is Lipschitz continuous for either fixed $x$ or $y$ with constant $M$. This implies that the norms of the partial gradients are bounded by $M$:
\[
\|\nabla_xf(x,y)\|\le M,\quad \|\nabla_yf(x,y)\|\le M.
\]

In particular, we consider the settings when the partial gradients $\nabla_xf(x,y)$ and $\nabla_yf(x,y)$ (and any higher-order information) are not available. Furthermore, the {\it exact} evaluation of the function value itself is also not available; instead, we can only access a {\it stochastic oracle} $\hat{f}(x,y,\xi)$, which satisfies the following assumption:
\begin{eqnarray}
&&\left\{ \begin{array}{l}
\mathbb{E}\left[\hat{f}(x,y,\xi)\right]=f(x,y),\\ \\
\mathbb{E}\left[\nabla_x\hat{f}(x,y,\xi)\right]=\nabla_xf(x,y),\\ \\
\mathbb{E}\left[\nabla_y\hat{f}(x,y,\xi)\right]=\nabla_yf(x,y),\\ \\
\mathbb{E}\left[\|\nabla_x\hat{f}(x,y,
\xi)-\nabla_xf(x,y)\|^2\right]\le \sigma^2,\\ \\
\mathbb{E}\left[\|\nabla_y\hat{f}(x,y,
\xi)-\nabla_yf(x,y)\|^2\right]\le \sigma^2.
\end{array}\right.\label{so-x-mean}
\end{eqnarray}

Now, we shall use the so-called {\it smoothing}\/ technique to approximate the first-order information, which then enables us to apply the proposed stochastic methods for VI, which includes the saddle-point model as a special case. In particular, we use a randomized smoothing scheme using uniform distributions $U_b$/$V_b$ over the unit Euclidean ball $B$ in the $\mathbb{R}^n$/$\mathbb{R}^m$ space, respectively. The smoothing functions with parameters $\rho_x,\rho_y> 0$ are defined as follows:
\begin{eqnarray}
f_{\rho_x}(x,y) &:=& \mathbb{E}_{u\sim U_b}\left[f(x+\rho_xu,y)\right]=\frac{1}{\alpha(n)}\int_Bf(x+\rho_xu,y)du,\nonumber\\
f_{\rho_y}(x,y) &:=& \mathbb{E}_{v\sim V_b}\left[f(x,y+\rho_yv)\right]=\frac{1}{\alpha(m)}\int_Bf(x,y+\rho_yv)dv,\nonumber
\end{eqnarray}
where $\alpha(n)$/$\alpha(m)$ is the volume of the unit ball in $\mathbb{R}^n$/$\mathbb{R}^m$.

Let us summarize main properties of the smoothing functions $f_{\rho_x},\, f_{\rho_y}$ below:

\begin{lemma}
\label{lem:smooth-function}
Let $U_{S_p}/V_{S_p}$ be the uniform distribution on the unit sphere $S_p$ in $\mathbb{R}^n/\mathbb{R}^m$. Then,
\begin{enumerate}
    \item The smoothing functions $f_{\rho_x},f_{\rho_y}$ are continuously differentiable, and their partial gradients $\nabla_xf_{\rho_x},\nabla_yf_{\rho_y}$ can be expressed as:
    \begin{eqnarray}
    \nabla_xf_{\rho_x}(x,y)&:=& \mathbb{E}_{u\sim U_{S_p}}\left[\frac{n}{\rho_x}f(x+\rho_xu,y)u\right]=\frac{1}{\beta(n)}\int_{u\in S_p}\frac{n}{\rho_x}\left(f(x+\rho_xu,y)-f(x,y)\right)udu,\nonumber\\
    \nabla_yf_{\rho_y}(x,y)&:=&\mathbb{E}_{v\sim V_{S_p}}\left[\frac{m}{\rho_y}f(x,y+\rho_yv)v\right]=\frac{1}{\beta(m)}\int_{v\in S_p}\frac{m}{\rho_y}\left(f(x,y+\rho_yv)-f(x,y)\right)vdv,\nonumber
    \end{eqnarray}
    where $\beta(n)/\beta(m)$ is the surface area of the unit sphere in $\mathbb{R}^n/\mathbb{R}^m$.
    %Furthermore, the partial gradients are Lipschitz continuous with constants $L_{\rho_x},L_{\rho_y}$ for fixed $y/x$. Denote $L_{\rho}=\max\{L_{\rho_x},L_{\rho_y}\}$, then we have $L_{\rho}\le L$.
    \item For any $x\in\mathbb{R}^n$ and any $y\in\mathbb{R}^m$, we have:
    \begin{eqnarray}
    %&& |f_{\rho_x}(x,y)-f(x,y)|\le \frac{L\rho_x^2}{2},\nonumber\\
    %&&|f_{\rho_y}(x,y)-f(x,y)|\le \frac{L\rho_y^2}{2},\nonumber\\
    &&\left\{\begin{array}{l}
         \|\nabla_xf_{\rho_x}(x,y)-\nabla_xf(x,y)\|\le\frac{\rho_xnL}{2},\\ \\
        \|\nabla_yf_{\rho_y}(x,y)-\nabla_yf(x,y)\|\le \frac{\rho_ymL}{2},%\label{smooth-grad-bd-y}\\
    \end{array}\right.\label{smooth-grad-bd}\\ \nonumber \\
    &&\left\{
    \begin{array}{l}
        \mathbb{E}_{u\sim U_{S_p}}\left[\left\|\frac{n}{\rho_x}\left(f(x+\rho_xu,y)-f(x,y)\right)u\right\|^2\right]\le 2n\|\nabla_xf(x,y)\|^2+\frac{\rho_x^2L^2n^2}{2},\\
        \\
        \mathbb{E}_{v\sim V_{S_p}}\left[\left\|\frac{m}{\rho_y}\left(f(x,y+\rho_yv)-f(x,y)\right)v\right\|^2\right]\le 2m\|\nabla_yf(x,y)\|^2+\frac{\rho_y^2L^2m^2}{2}.
    \end{array}
    \right.\label{smooth-x-sm}
    \end{eqnarray}
\end{enumerate}
\end{lemma}
\begin{proof}
For Statement 1, cf.~\cite{shalev2011online} (Lemma 4.4), and for Statement 2, cf.~\cite{gao2018low} (proofs for Propositions 2.7.5 and 2.7.6). Note that the proofs for the minimax function follows simply by fixing one of the two variables.
\end{proof}

We are now ready to define the {\it stochastic zeroth-order gradient}\/ as follows:
\begin{eqnarray}
&&\left\{\begin{array}{l}
    F_{\rho_x}(x,y,\xi,u) := \frac{n}{\rho_x}\left(\hat{f}(x+\rho_xu,y,\xi)-\hat{f}(x,y,\xi)\right)u,\\
    F_{\rho_y}(x,y,\xi,v) := \frac{m}{\rho_y}\left(\hat{f}(x,y+\rho_yv,\xi)-\hat{f}(x,y,\xi)\right)v,
\end{array}\right.\label{sto-0-grad}
\end{eqnarray}
where $u$ and $v$ are the uniformly distributed random vectors over the unit spheres in $\mathbb{R}^n$ and $\mathbb{R}^m$ respectively.

The next lemma shows that such stochastic zeroth-order gradients are unbiased with respect to the gradients of the smoothing functions and have uniformly bounded variance.

\begin{lemma}
\label{lem:sto-0-grad}
The stochastic zeroth-order gradients defined in \eqref{sto-0-grad} are unbiased and have bounded variance for all $(x,y)\in\mathcal{X}\times\mathcal{Y}$:
\begin{eqnarray}
\left\{\begin{array}{l}
    \mathbb{E}_{\xi,u}\left[F_{\rho_x}(x,y,\xi,u)\right]=\nabla_xf_{\rho_x}(x,y),\\ \\
    \mathbb{E}_{\xi,v}\left[F_{\rho_y}(x,y,\xi,v)\right]=\nabla_yf_{\rho_y}(x,y),
\end{array}\right.\label{SZG-mean}
\end{eqnarray}
and
\begin{eqnarray}
\left\{\begin{array}{l}
     \mathbb{E}_{\xi,u}\left[\|F_{\rho_x}(x,y,\xi,u)-\nabla_xf_{\rho_{x}}(x,y)\|^2\right]\le \tilde{\sigma}^2,\\
     \\
     \mathbb{E}_{\xi,v}\left[\|F_{\rho_y}(x,y,\xi,v)-\nabla_yf_{\rho_{y}}(x,y)\|^2\right]\le \tilde{\sigma}^2,
\end{array}\right.
\label{SZG-var}
\end{eqnarray}
where
\[
\tilde{\sigma}^2 = 2\cdot\max\left\{nM^2+n\sigma^2+n^2\rho_x^2L^2,\, mM^2+m\sigma^2+m^2\rho_y^2L^2 \right\}.
\]
\end{lemma}
\begin{proof}
See Appendix \ref{proof:sto-0-grad}.
\end{proof}

Before applying the stochastic extra-point/extra-momentum scheme to solve~\eqref{spp}, let us first introduce the connections between these two models. %define some notations to ease the burden in the following analysis.
%Since we are applying methods originally for VIs,
As we regard the saddle-point model as a special case of VI,
we shall treat the variables $x,y$ in the saddle-point problem as one variable and denote $\mathcal{Z}=\mathcal{X}\times\mathcal{Y}, z=(x,y)$. Additionally, we define:
\begin{eqnarray*}
&& F(z) := \begin{pmatrix}
\nabla_xf(x,y)\\
-\nabla_yf(x,y)
\end{pmatrix},
% \quad \hat{F}(z,\xi)=\begin{pmatrix}
% \nabla_x\hat{f}(x,y,\xi)\\
% -\nabla_y\hat{f}(x,y,\xi)
% \end{pmatrix},\\
\quad F_{\rho}(z) := \begin{pmatrix}
\nabla_xf_{\rho_x}(x,y)\\
-\nabla_yf_{\rho_y}(x,y)
\end{pmatrix},\quad \hat{F}_{\rho}(z,u,v,\xi) := \begin{pmatrix}
F_{\rho_x}(x,y,\xi,u)\\
-F_{\rho_y}(x,y,\xi,v)
\end{pmatrix}.
\end{eqnarray*}
These terms correspond to the gradient of $f(x,y)$, the gradient of the smoothing functions $f_{\rho_x}(x,y)$ and $f_{\rho_{y}}(x,y)$, and the stochastic zeroth-order gradient, respectively. Note that we have flipped the sign on partial gradient correspond to $y$ to account for the {\it concavity}\/ of $f$ with respect to $y$.

Finally, as we shall use a sample size of $t_k\in \mathbb{N}$ (a natural number) at iteration $k$, we reserve the subscripts for the random vectors $\xi,u,v$ for the sample index $(i),\, i=1,...,t_k$, and denote:
\[
\hat{F}^k_{\rho}(z^k)=\frac{1}{t_k}\sum\limits_{i=1}^{t_k}\hat{F}_{\rho}(z^k,u^k_{(i)},v^k_{(i)},\xi^k_{(i)}).
\]
In the above definition we suppress the notation of the random vectors $u,v,\xi$ on the LHS for cleaner presentation. Note that by the {\it law of large numbers}, together with \eqref{SZG-mean}-\eqref{SZG-var}, we have
\begin{eqnarray}
&& \mathbb{E}\left[\hat{F}^k_{\rho}(z^k)\right]=F_{\rho}(z^k),\label{batch-mean}\\
&& \mathbb{E}\left[\|\hat{F}^k_{\rho}(z^k)-F_{\rho}(z^k)\|^2\right]\le \frac{2\tilde{\sigma}^2}{t_k}.\label{batch-var}
\end{eqnarray}

\subsection{Sample complexity analysis: stochastic zeroth-order extra-point method}
\label{sec:sample-complex}
Recall that our objective is
\[
\min\limits_{x\in\mathcal{X}}\max\limits_{y\in\mathcal{Y}}f(x,y).
\]
With only noisy function value $\hat{f}(x,y,\xi)$ accessible, we propose the {\it stochastic zeroth-order extra-point method} that updates $(x,y):=z$ simultaneously with the following update rule:
\begin{eqnarray}
&&\left\{
\begin{array}{lcl}
     z^{k+0.5} &:=& P_{\mathcal{Z}}\left(z^k+\beta(z^k-z^{k-1})-\eta \hat{F}^k_{\rho}(z^k)\right), \\
     z^{k+1} &:=& P_{\mathcal{Z}}\left(z^k-\alpha \hat{F}^{k+0.5}_{\rho}(z^{k+0.5})+\gamma(z^k-z^{k-1})-\tau\left(\hat{F}^k_{\rho}(z^{k})-\hat{F}^{k-1}_{\rho}(z^{k-1})\right)\right). \\
\end{array}
\right.\label{szo-exp-update}
\end{eqnarray}

Compare the above update with its original variant in \eqref{sto-exp-update} for solving stochastic VI, the update direction $\hat{F}(z^k)$ is replaced by the averaged stochastic zeroth-order gradient $\hat{F}^k_{\rho}(z^k)$ with sample size $t_k$ (similarly $\hat{F}(z^{k+0.5})$ is replaced by $\hat{F}^{k+0.5}_{\rho}(z^{k+0.5})$ with sample size $t_{k+0.5}$). This circumvents the inaccessible first-order information and equips us with appropriate tools to reduce the variance and achieve overall linear convergence.

The next lemma establishes the relational inequality between the subsequent iterates in terms of the expected distance to the solution $d_k=\mathbb{E}\left[\|z^k-z^*\|^2\right]$, similar to what we did in Section \ref{sec:sto-exp}. The differences lie in the corresponding stochastic error terms shown below. Note that in each iteration we take two batches of samples $t_k$ and $t_{k+0.5}$. The batch size $t_{k-1}$ also appears because the iterate $z^{k-1}$ is used in each iteration.
\begin{lemma}
\label{lem:sto-0-exp-relation}
For the sequences $\{z^k\mid k=0,1,...\}$ and $\{z^{k+0.5}\mid k=0,1,...\}$ generated from the stochastic zeroth-order extra-point method \eqref{szo-exp-update}, the following inequality holds:
\begin{eqnarray}
&&(1-4|\gamma-\beta|-\tau L)d_{k+1}\nonumber\\
&\le& \left(1+4\gamma+6|\gamma-\beta|+4\tau L-\alpha\mu\right)d_k\nonumber\\
&&+\left(2|\gamma-\beta|+2\gamma+4\tau L\right)d_{k-1}\nonumber\\
&&+\left(\alpha^2L^2+2|\gamma-\beta|+2\gamma+2\alpha\mu-1\right)\mathbb{E}\left[\|z^{k+0.5}-z^k\|^2\right]\nonumber\\
&&+16\tilde{\sigma}^2\left(\left(\alpha^2+\frac{\tau}{L}\right)\frac{1}{t_k}+\frac{\alpha^2}{t_{k+0.5}}+\frac{\tau}{Lt_{k-1}}\right)+\alpha LD\sqrt{\rho_x^2n^2+\rho_y^2m^2}\nonumber\\
&&+ 4(\tau L+\alpha^2L^2)(\rho_x^2n^2+\rho_y^2m^2)\nonumber\\
&&+2(\eta-\alpha)\mathbb{E}\left[(z^{k+1}-z^{k+0.5})^\top\hat{F}^k_{\rho}(z^k)\right].\label{szo-exp-relation}
\end{eqnarray}
\end{lemma}
\begin{proof}
See Appendix \ref{proof:sto-0-exp-relation}.
\end{proof}

With the relational inequality in Lemma \ref{lem:sto-0-exp-relation}, we shall adopt the {\it same} conditions: \eqref{sto-exp-cond-1}, \eqref{sto-exp-t}, and \eqref{sto-exp-cond-t} for the parameters $\alpha,\beta,\gamma,\eta,\tau$. Therefore, the results in Theorem \ref{th:sto-exp-result} are directly applicable. In addition, we are now equipped with the variable sample size $t_k$/$t_{k+0.5}$ to control the variance terms, as well as the smoothing parameters $\rho_x,\rho_y$ to control the bias terms
\[
\alpha LD\sqrt{\rho_x^2n^2+\rho_y^2m^2}+4(\tau L+\alpha^2L^2)(\rho_x^2n^2+\rho_y^2m^2).
\]

We shall utilize the example in Proposition \ref{prop:sto-exp-example} to analyze the sample complexity of the proposed method. The result is provided in the next proposition:
\begin{proposition}[Sample complexity result 1]
\label{prop:szo-exp-sample}
The stochastic zeroth-order extra-point method \eqref{szo-exp-update} with the following parameter choice:
\begin{eqnarray}
(\alpha,\beta,\gamma,\eta,\tau)=\left(\frac{1}{4L},\frac{1}{64\kappa},\frac{1}{64\kappa},\frac{1}{4L},\frac{1}{64L\kappa}\right)\nonumber
\end{eqnarray}
and
\[
t_{k+0.5} = t_k = \ceil[\Bigg]{K\left(1-\frac{1}{256\kappa}\right)^{-(k+1)}},\quad \rho_x=\frac{1}{\sqrt{2}n\kappa}\left(1-\frac{1}{256\kappa}\right)^{K},\quad \rho_y=\frac{1}{\sqrt{2}m\kappa}\left(1-\frac{1}{256\kappa}\right)^{K},
\]
where $K$ is the iteration count decided in advance, outputs $z^K$ such that $d_K=\mathbb{E}\left[\|z^K-z^*\|^2\right]\le \epsilon$, with $K=\mathcal{O}\left(\kappa\ln\left(\frac{1}{\epsilon}\right)\right)$, and the total sample complexity of the procedure is
\[
\sum\limits_{k=0}^{K-1}(t_k+t_{k+0.5})=\mathcal{O}\left(\frac{\kappa^2}{\epsilon}\ln\left(\frac{1}{\epsilon}\right)\right).
\]
\end{proposition}
\begin{proof}
See Appendix \ref{proof:prop-szo-exp-sample}.
\end{proof}

\subsection{Sample complexity: stochastic zeroth-order extra-momentum method}

%Following the similar logic,
Next we consider the {\it stochastic zeroth-order extra-momentum method}, with one projection per each iteration:
\begin{eqnarray}
z^{k+1} := P_{\mathcal{Z}}\left(z^k-\alpha \hat{F}^k_{\rho}(z^k)+\gamma(z^k-z^{k-1})-\tau\left(\hat{F}^k_{\rho}(z^k)-\hat{F}^{k-1}_{\rho}(z^{k-1})\right)\right). \label{szo-exm-update}
\end{eqnarray}

The relational inequality, similar to Lemma \ref{lem:sto-exm-relation}, is established in the next lemma:
\begin{lemma}
\label{lem:szo-exm-relation}
For the sequence $\{z^k\mid k=0,1,...\}$ generated from the stochastic zeroth-order extra-momentum method \eqref{szo-exm-update}, the following inequality holds:
\begin{eqnarray}
&&\mathbb{E}\left[\left(\frac{1}{2}+\frac{\alpha\mu}{2}-\frac{\gamma}{2}\right)\|z^{k+1}-z^*\|^2+\alpha(z^{k+1}-z^*)^\top\left(\hat{F}^k_{\rho}(z^k)-\hat{F}^{k+1}_{\rho}(z^{k+1})\right)+\frac{1}{4}\|z^{k+1}-z^k\|^2\right]\nonumber\\
&\le& \mathbb{E}\left[\frac{1}{2}\|z^k-z^*\|^2+\tau(z^k-z^*)^\top\left(\hat{F}^{k-1}_{\rho}(z^{k-1})-\hat{F}^k_{\rho}(z^{k})\right)+\left(2\tau^2L^2+\frac{\gamma}{2}\right)\|z^k-z^{k-1}\|^2\right]\nonumber\\
&&+16\tau^2\sigma^2\left(\frac{1}{t_k}+\frac{1}{t_{k-1}}\right)+L^2\left(4\tau^2+\frac{\alpha}{8\mu}\right)(\rho_x^2n^2+\rho_y^2m^2).\label{szo-exm-relation}
\end{eqnarray}
\end{lemma}
\begin{proof}
See Appendix \ref{proof:szo-exm-relation}.
\end{proof}

With the {\it same} condition as in \eqref{sto-exm-condition-par} for the parameters $\alpha,\gamma,\tau$, we can derive the similar bound to \eqref{sto-exm-relation-2} (with $\hat{F}(z^k)$ replaced with $\hat{F}^k_{\rho}(z^k)$ and with the new stochastic error expression) and define the potential function:
\[
V_k = \mathbb{E}\left[\frac{1}{2}\|z^k-z^*\|^2+\tau(z^k-z^*)^\top\left(\hat{F}^{k-1}_{\rho}(z^{k-1})-\hat{F}^k_{\rho}(z^{k})\right)+\left(2\tau^2L^2+\frac{\gamma}{2}\right)\|z^k-z^{k-1}\|^2\right].
\]

Therefore, the following inequality holds:
\begin{eqnarray}
\left(1+\frac{\theta}{\kappa}\right)V_{k+1}\le V_k+16\tau^2\sigma^2\left(\frac{1}{t_k}+\frac{1}{t_{k-1}}\right)+L^2\left(4\tau^2+\frac{\alpha}{8\mu}\right)(\rho_x^2n^2+\rho_y^2m^2),\label{szo-exm-potential-relation}
\end{eqnarray}
and we can apply the results directly from Theorem \ref{th:sto-exm-result}. In addition, with increasing sample sizes $t_k$ and the smoothing parameters $\rho_x,\rho_y$, we are able to control the bias and the variance terms in the above inequality. We give the results of sample complexity in the next proposition.

\begin{proposition}[Sample complexity result 2]
\label{prop:szo-exm-sample}
The stochastic zeroth-order extra-momentum method \eqref{szo-exm-update} with the following parameter choice:
\[
(\alpha,\tau,\gamma)=\left(\frac{1}{4L},\frac{\alpha}{1+\frac{\theta}{\kappa}},\frac{1}{8(\kappa+\theta)}\right),\quad \theta=\frac{1}{8},
\]
and
\[
t_k=\ceil[\Bigg]{K\left(1-\frac{1}{8\kappa+1}\right)^{-k}},\quad \rho_x=\frac{1}{\sqrt{2}n\kappa}\left(1-\frac{1}{8\kappa+1}\right)^{\frac{K}{2}},\quad \rho_y=\frac{1}{\sqrt{2}m\kappa}\left(1-\frac{1}{8\kappa+1}\right)^{\frac{K}{2}},
\]
where $K$ is the iteration count decided in advance, outputs $z^K$ such that $d_K=\mathbb{E}\left[\|z^K-z^*\|^2\right]\le \epsilon$, with $K=\mathcal{O}\left(\kappa\ln\left(\frac{1}{\epsilon}\right)\right)$ and the total sample complexity of the procedure is
\[
\sum\limits_{k=0}^{K-1}t_k=\mathcal{O}\left(\frac{\kappa^2}{\epsilon}\ln\left(\frac{1}{\epsilon}\right)\right).
\]
\end{proposition}
\begin{proof}
See Appendix \ref{proof:szo-exm-sample}.
\end{proof}

Remark that the sample complexity of VS-Ave proposed in \cite{jalilzadeh2019proximal} for solving stochastic strongly monotone VI is given by $\mathcal{O}\left(\left(\frac{\kappa^2}{\epsilon}\right)^{\beta}\right)$ for some $\beta>1$, and its limiting case ($\beta\rightarrow 1$) differs from our results given in Proposition \ref{prop:szo-exp-sample} and \ref{prop:szo-exm-sample} only by a factor of $\mathcal{O}\left(\ln\left(\frac{1}{\epsilon}\right)\right)$.

\section{Numerical Experiments}
\label{sec:numerical}

In this section, we conduct an experiment that models a regularized two-player zero-sum matrix game with some {\it uncertain}\/ payoff matrix. In particular, the payoff matrix $A_\xi$ is randomly distributed and can only be sampled for each (mixed) strategy. The problem can be formulated as follows:
\begin{eqnarray}
\min\limits_{x\in\mathbb{R}^n}\max\limits_{y\in\mathbb{R}^m}&&f(x,y)= \mathbb{E} \left[ \frac{\lambda_x}{2}\|x\|^2+x^\top A_\xi y-\frac{\lambda_y}{2}\|y\|^2 \right] \label{numerical-prob}\\
\mbox{s.t.}&& \sum\limits_{i=1}^nx_i=1,\quad\sum\limits_{j=1}^my_j=1,\nonumber\\
&&x,y\ge0.\nonumber
\end{eqnarray}

The experiment consists of two parts. In the first part, the random matrix $A_{\xi}$ is sampled element-wise from i.i.d.\ normal distribution, $A_{\xi}\sim N \left(A_0,\sigma^2I_{(n+m)}\right)$, where $A_0$ is pre-determined and randomly generated as follows. We partition $A_0=\begin{pmatrix}A_{11}&A_{12}\\A_{21}&A_{22}\end{pmatrix}$, with each block matrix $A_{ij}\in\mathbb{R}^{\frac{n}{2}\times\frac{m}{2}}$. Each entry in $A_{ij}$ is generated randomly from Unif$(a_{ij}-b_{ij},a_{ij}+b_{ij})$, where $a_{ij}$ and $b_{ij}$ are randomly generated from Unif$(-30,30)$ and Unif$(0,30)$ respectively. The problem parameters are set as: $n=10$, $m=20$, $\sigma^2=0.5$, $\lambda_x=\lambda_y=1$, $\kappa=\frac{L}{\mu}\approx161$, where $L$ and $\mu$ are the largest and the smallest singular values of the Jacobian matrix $\begin{pmatrix}\lambda I_n & A_0\\-A_0^\top&\lambda I_m\end{pmatrix}$ respectively. The smoothing parameters $\rho_x,\rho_y$ are set to be in the order $10^{-8}$, the total iteration $K$ is set to 1485, and the sample size $t_k=k$ at iteration $k$.

In the second part, the random matrix $A_{\xi}$ is sampled element-wise from i.i.d.\ {\it log-normal} distribution, which is known as to have a {\it fat-tail}. It is used to model multiplicative random variables that take positive values. % and has applications in engineering, medicine, economics, to name a few.
We reuse the parameters $A_0$ and $\sigma^2$ from the first part. In particular, the samples are generated by $A_{\xi}=e^{\frac{A_0}{10}+\sigma Z}$, where $Z$ is sampled element-wise from i.i.d.\ standard normal distribution $N(0,1)$. Therefore, the mean of such distribution is given by $A_0'=\mathbb{E}\left[A_{\xi}\right]=e^{\frac{A_0}{10}+\frac{\sigma^2}{2}}$. We have $\kappa\approx 146$ and set $K=1345$ in this part, and other parameters remain the same as the first part ($\rho_x,\rho_y$ are in the same order).

In both parts of the experiment, we first solve the deterministic problem with the mean payoff matrix $A_0$ ($A_0'$ for the second part) and denote the solution as $(x^*,y^*)$. We then implement the two proposed methods: the stochastic zeroth-order extra-point method and the stochastic zeroth-order extra-momentum method. In addition, we compare these two methods with other first order methods: the extra-gradient method, the OGDA method, and the VS-Ave method (proposed in \cite{jalilzadeh2019proximal}, which is a variance-reduced stochastic extension of Nesterov's method \cite{nesterov2006solving}), all equipped with the same stochastic zeroth-order oracle. The results are shown in the following two figures, where the left plot shows the result from one experiment and the right shows the result from average of ten experiments. The parameters for each algorithm are manually tuned except for VS-Ave, where we adopt the recursive rule as proposed in its original paper. The results show that the two newly proposed methods exhibit comparable (or slightly improved) performance to the stochastic extra-gradient/OGDA method in this particular example of application.

\begin{figure}[htbp]
\centering
\begin{minipage}[t]{0.48\textwidth}
\centering
\includegraphics[width=7cm]{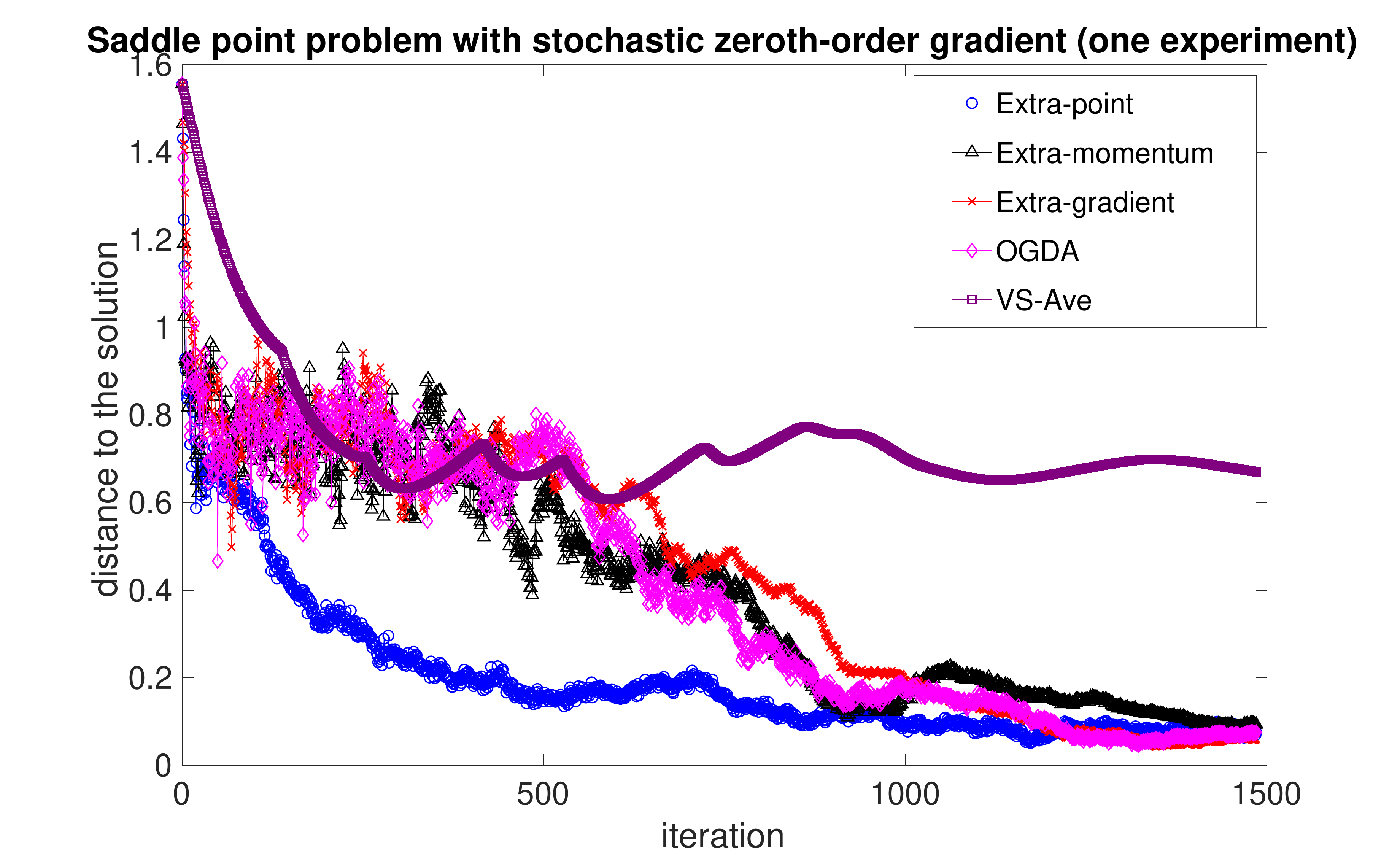}
%\label{fig:epmt-1}
\end{minipage}
\begin{minipage}[t]{0.48\textwidth}
\centering
\includegraphics[width=7cm]{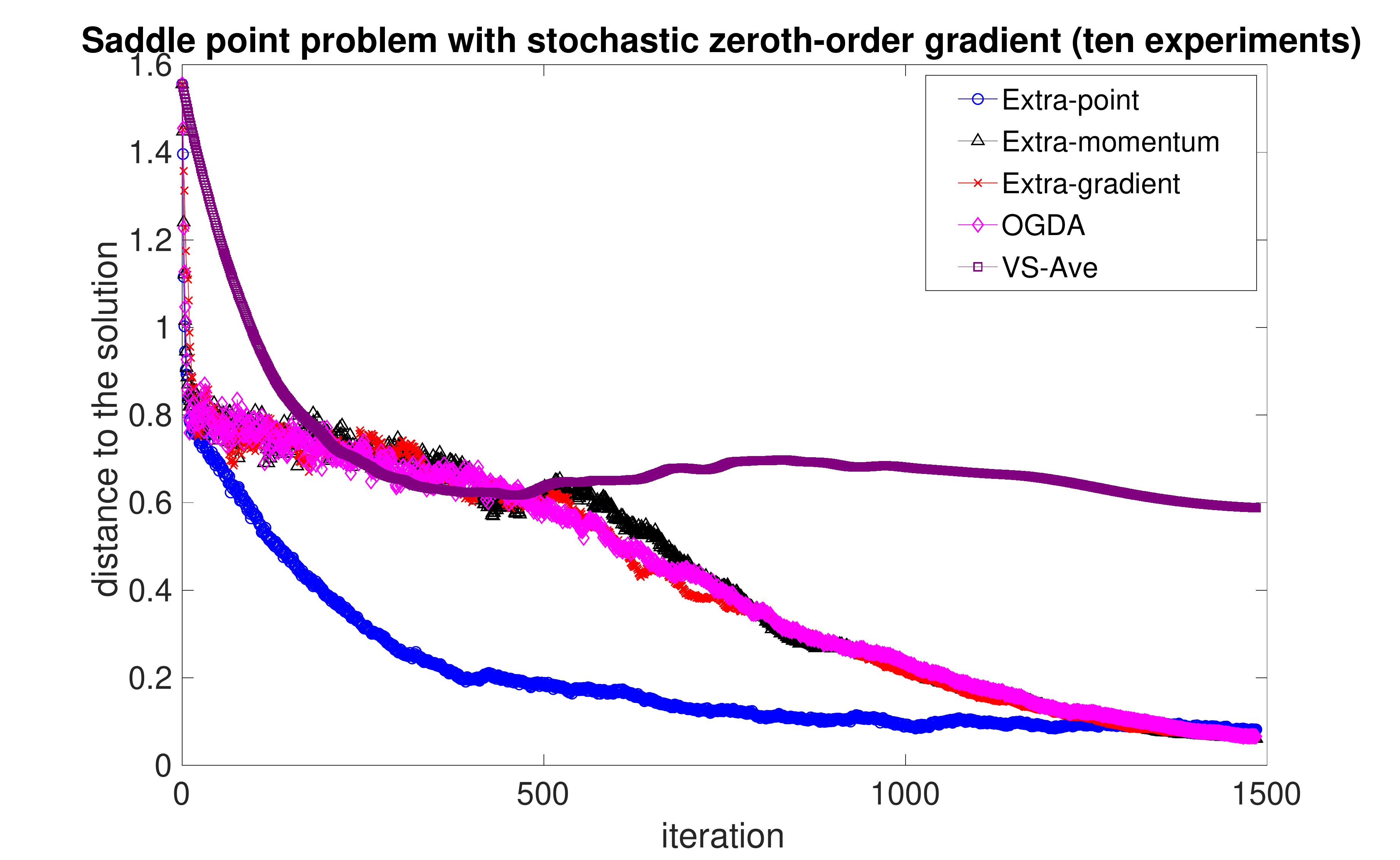}
%\label{fig:epmt-10}
\end{minipage}
\caption{Normal distributed payoff matrix $A$}
\label{fig:normal}
\end{figure}

\begin{figure}[htbp]
\centering
\begin{minipage}[t]{0.48\textwidth}
\centering
\includegraphics[width=7cm]{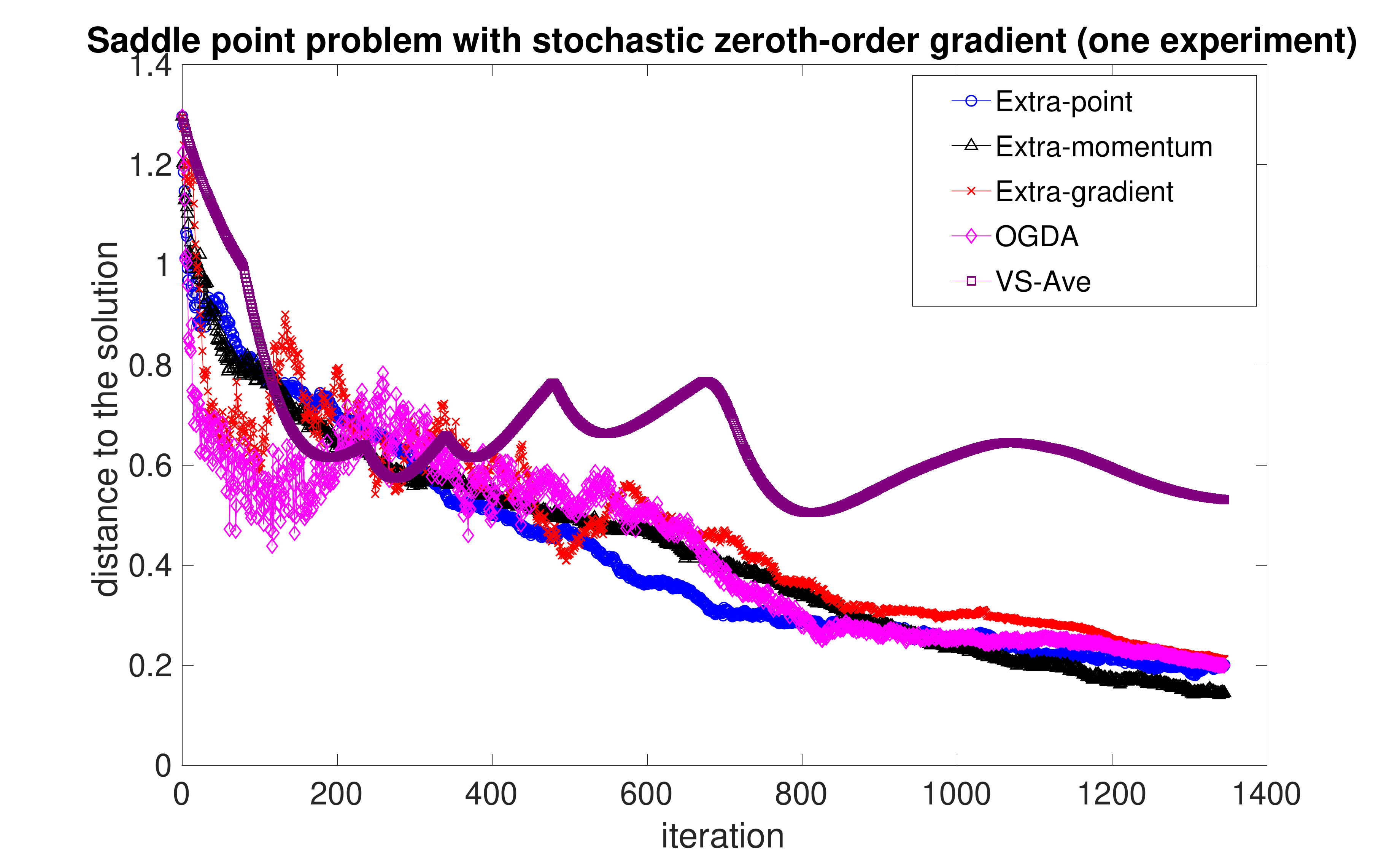}
%\label{fig:epmt-1-log-normal}
\end{minipage}
\begin{minipage}[t]{0.48\textwidth}
\centering
\includegraphics[width=7cm]{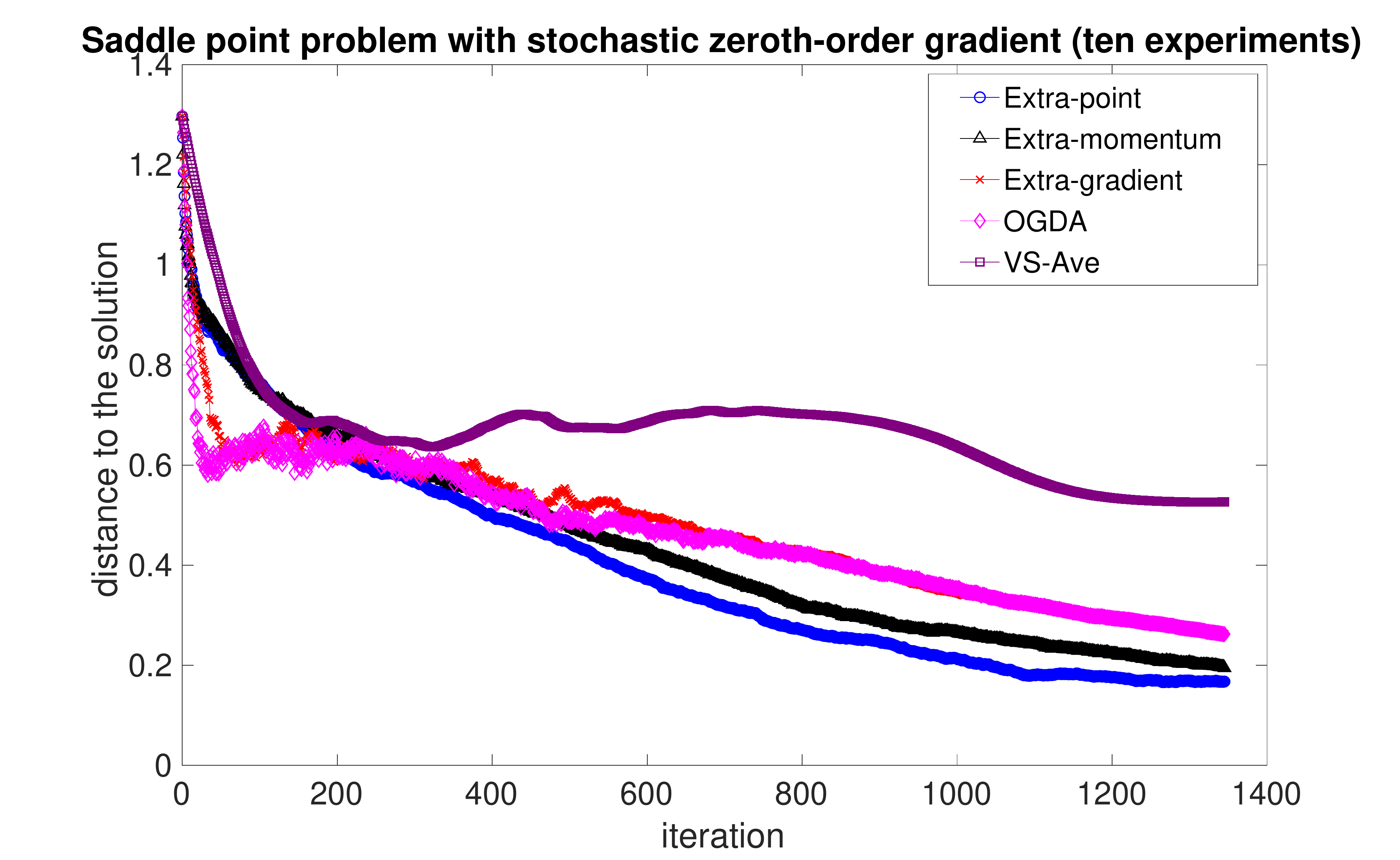}
%\label{fig:epmt-10-log-normal}
\end{minipage}
\caption{Log-normal distributed payoff matrix $A$}
\label{fig:log-normal}
\end{figure}

\section{Conclusions}
\label{sec:conclusion}
This paper proposes two new schemes of stochastic first-order methods to solve strongly monotone VI problems: the stochastic extra-point scheme and the stochastic extra-momentum scheme. The first scheme features a high flexibility in the configuration of parameter choices that can be tailored to different problem classes.
%With a tradeoff in generalization (by specifically combining two of the search directions from the first method),
The second scheme is less general in the choice of search directions. However, it has the advantage of maintaining a single sequence throughout the iterations. Therefore, it requires only one projection per iteration, as opposed to most other first method that maintains an extra iterative sequence. Both methods achieve optimal iteration complexity bound, provided that the stochastic gradient oracle allows the variance to be controllable. The application of these two schemes to solve stochastic black-box saddle-point problem is also presented. Through a randomized smoothing scheme, the stochastic oracles required in these two schemes can be constructed via stochastic zeroth-order gradient approximation. The variance is thus controllable by mini-batch sampling with linearly increasing sample sizes per iteration, and the sample complexity results are derived. Preliminary numerical experiments show an improved (or at least comparable) performance of the proposed schemes to other existing methods.

\printbibliography
% \bibliographystyle{acm}
% \bibliography{plan}

\begin{appendices}
\section{Proof of technical lemmas and theorems}
\label{appendix:proofs}
\subsection{Proof of Lemma \ref{lem:sto-exp-relation}}

First of all, by the 1-co-coerciveness (cf.\ e.g.\ Proposition 4.4 in \cite{bauschke2011convex}) of the projection operator $P_{\mathcal{Z}}$, we have
\begin{eqnarray}
\|z^{k+1}-z^*\|^2 &=& \|P_{\mathcal{Z}}\left(z^k-\alpha\hat{F}(z^{k+0.5})+\gamma(z^k-z^{k-1})-\tau\left(\hat{F}(z^k)-\hat{F}(z^{k-1})\right)\right)-P_{\mathcal{Z}}\left(z^*\right)\|^2 \nonumber\\
&\le& (z^{k+1}-z^*)^\top\left(z^k-\alpha\hat{F}(z^{k+0.5})+\gamma(z^k-z^{k-1})-\tau\left(\hat{F}(z^k)-\hat{F}(z^{k-1})\right)-z^*\right) \nonumber\\
&=& \frac{1}{2}\|z^{k+1}-z^*\|^2+\frac{1}{2}\|z^k-z^*\|^2-\frac{1}{2}\|z^{k+1}-z^k\|^2-\tau(z^{k+1}-z^*)^\top\left(\hat{F}(z^k)-\hat{F}(z^{k-1})\right)\nonumber\\
&&+(z^{k+1}-z^*)^\top\left(-\alpha\hat{F}(z^{k+0.5})+\gamma(z^k-z^{k-1})\right).\label{bd-first-step}
\end{eqnarray}

We shall first decompose the last term in the above inequality as
\begin{eqnarray}
&& (z^{k+1}-z^*)^\top\left(-\alpha\hat{F}(z^{k+0.5})+\gamma(z^k-z^{k-1})\right) \nonumber\\
&=& (z^{k+1}-z^{k+0.5}+z^{k+0.5}-z^*)^\top\left(-\alpha\hat{F}(z^{k+0.5})+\gamma(z^k-z^{k-1})\right)\nonumber\\
&=& \underbrace{(z^{k+1}-z^{k+0.5})^\top\left(-\eta\hat{F}(z^k)+\beta(z^k-z^{k-1})\right)}_{(a)}\nonumber\\
&& + \underbrace{(z^{k+1}-z^{k+0.5})^\top\left(-\alpha\hat{F}(z^{k+0.5})+\eta\hat{F}(z^k)+(\gamma-\beta)(z^k-z^{k-1})\right)}_{(b)}\nonumber\\
&&+\underbrace{(z^{k+0.5}-z^*)^\top\left(-\alpha\hat{F}(z^{k+0.5})+\gamma(z^k-z^{k-1})\right)}_{(c)} . \label{bd-last}
\end{eqnarray}

Let us first use the optimality condition of $z^{k+0.5}$ to bound term $(a)$:
\[
\langle z^{k+0.5}-z^k-\beta(z^k-z^{k-1})+\eta\hat{F}(z^k),z-z^{k+0.5}\rangle\geq 0,\quad\forall z\in\mathcal{Z}.
\]
Taking $z=z^{k+1}$, we get
\begin{eqnarray}
(a)&\le& \frac{1}{2}\|z^{k+1}-z^k\|^2-\frac{1}{2}\|z^{k+0.5}-z^k\|^2-\frac{1}{2}
\|z^{k+1}-z^{k+0.5}\|^2. \label{bd-last-a}
\end{eqnarray}

We can also establish the bound for $(b)$:
\begin{eqnarray}
(b)&=& (z^{k+1}-z^{k+0.5})^\top\left(-\alpha\hat{F}(z^{k+0.5})+\alpha \hat{F}(z^k)-\alpha \hat{F}(z^k)+\eta\hat{F}(z^k)+(\gamma-\beta)(z^k-z^{k-1})\right)\nonumber\\
&\le& \alpha\|z^{k+1}-z^{k+0.5}\|\|\hat{F}(z^k)-\hat{F}(z^{k+0.5})\|+(\eta-\alpha)(z^{k+1}-z^{k+0.5})^\top\hat{F}(z^k)\nonumber\\
&&+(\gamma-\beta)(z^{k+1}-z^{k+0.5})^\top(z^k-z^{k-1}). \nonumber
\end{eqnarray}

Note the following bound from the Lipschitz continuity:
\begin{eqnarray}
\|\hat{F}(z)-\hat{F}(z')\|&=& \|\hat{F}(z)-F(z)+F(z')-\hat{F}(z')+F(z)-F(z')\| \nonumber\\
&\le& \varepsilon_z+\varepsilon_{z'}+L\|z-z'\|,\label{sto-Lip-bd}
\end{eqnarray}
for any $z,z'\in\mathcal{Z}$, where we used the definition of the stochastic error term
\begin{eqnarray}
\varepsilon_z = \left\|F(z)-\hat{F}(z)\right\|.
\end{eqnarray}
Therefore,
\begin{eqnarray}
&&\alpha\|z^{k+1}-z^{k+0.5}\|\|\hat{F}(z^k)-\hat{F}(z^{k+0.5})\|\nonumber\\
&\le& \frac{1}{2}\left(\|z^{k+1}-z^{k+0.5}\|^2+\alpha^2\|\hat{F}(z^k)-\hat{F}(z^{k+0.5})\|^2\right)\nonumber\\
&\le& \frac{1}{2}\|z^{k+1}-z^{k+0.5}\|^2+\alpha^2(\varepsilon_{z^k}+\varepsilon_{z^{k+0.5}})^2+\alpha^2L^2\|z^k-z^{k+0.5}\|^2.\label{sto-exp-lip-bd-1}
\end{eqnarray}
Furthermore,
\begin{eqnarray}
&&(\gamma-\beta)(z^{k+1}-z^{k+0.5})^\top(z^k-z^{k-1})\nonumber\\
&\le& \frac{1}{2}|\gamma-\beta|\left(\|z^{k+1}-z^{k+0.5}\|^2+\|z^k-z^{k-1}\|^2\right)\nonumber\\
&\le& |\gamma-\beta|\left(\|z^{k+1}-z^{k}\|^2+\|z^{k+0.5}-z^k\|^2+\|z^k-z^{*}\|^2+\|z^{k-1}-z^*\|^2\right)\nonumber\\
&\le& |\gamma-\beta|\left(2\|z^{k+1}-z^{*}\|^2+2\|z^k-z^*\|^2+\|z^{k+0.5}-z^k\|^2+\|z^k-z^{*}\|^2+\|z^{k-1}-z^*\|^2\right)\nonumber\\
&=&|\gamma-\beta|\left(2\|z^{k+1}-z^{*}\|^2+3\|z^k-z^*\|^2+\|z^{k+0.5}-z^k\|^2+\|z^{k-1}-z^*\|^2\right). \nonumber
\end{eqnarray}

The resulting bound for $(b)$ becomes:
\begin{eqnarray}
(b)&\le& \frac{1}{2}\|z^{k+1}-z^{k+0.5}\|^2+\alpha^2(\epsilon_{z^k}+\epsilon_{z^{k+0.5}})^2+\left(\alpha^2L^2+|\gamma-\beta|\right)\|z^k-z^{k+0.5}\|^2\nonumber\\
&&+(\eta-\alpha)(z^{k+1}-z^{k+0.5})^\top\hat{F}(z^k)\nonumber\\
&&+|\gamma-\beta|\left(2\|z^{k+1}-z^{*}\|^2+3\|z^k-z^{*}\|^2+\|z^{k-1}-z^*\|^2\right)\label{bd-last-b}.
\end{eqnarray}

Next let us bound $(c)$ in \eqref{bd-last}. We have,
\begin{eqnarray}
(c)&=&-\alpha(z^{k+0.5}-z^*)^\top\hat{F}(z^{k+0.5})+\gamma(z^{k+0.5}-z^*)^\top(z^k-z^{k-1})\nonumber\\
&\le& -\alpha(z^{k+0.5}-z^*)^\top\hat{F}(z^{k+0.5})+\frac{1}{2}\gamma\left(\|z^{k+0.5}-z^*\|^2+\|z^k-z^{k-1}\|^2\right)\nonumber\\
&\le& -\alpha(z^{k+0.5}-z^*)^\top\hat{F}(z^{k+0.5})+\gamma\left(\|z^{k+0.5}-z^k\|^2+\|z^k-z^*\|^2+\|z^k-z^{*}\|^2+\|z^{k-1}-z^*\|^2\right)\nonumber\\
&=& -\alpha(z^{k+0.5}-z^*)^\top\hat{F}(z^{k+0.5})+\gamma\left(\|z^{k+0.5}-z^k\|^2+2\|z^k-z^*\|^2+\|z^{k-1}-z^*\|^2\right). \label{bd-last-c}
\end{eqnarray}

Combining the bounds for $(a),(b),(c)$ from \eqref{bd-last-a}, \eqref{bd-last-b}, and \eqref{bd-last-c},  it follows from \eqref{bd-last} that
\begin{eqnarray}
&& (z^{k+1}-z^*)^\top\left(-\alpha\hat{F}(z^{k+0.5})+\gamma(z^k-z^{k-1})\right) \nonumber\\
&\le& -\alpha(z^{k+0.5}-z^*)^\top\hat{F}(z^{k+0.5})+(\eta-\alpha)(z^{k+1}-z^{k+0.5})^\top\hat{F}(z^k)+\alpha^2(\varepsilon_{z^k}+\varepsilon_{z^{k+0.5}})^2\nonumber\\
&&+2\left(|\gamma-\beta|\right)\|z^{k+1}-z^*\|^2+\left(2\gamma+3|\gamma-\beta|\right)\|z^k-z^*\|^2+\left(|\gamma-\beta|+\gamma\right)\|z^{k-1}-z^*\|^2\nonumber\\
&&+\frac{1}{2}\|z^{k+1}-z^k\|^2+\left(\alpha^2L^2+|\gamma-\beta|+\gamma-\frac{1}{2}\right)\|z^{k+0.5}-z^k\|^2.\label{bd-last-final}
\end{eqnarray}

We also need to bound the following term in \eqref{bd-first-step}:
\begin{eqnarray}
&&-\tau(z^{k+1}-z^*)^\top\left(\hat{F}(z^k)-\hat{F}(z^{k-1})\right)\nonumber\\
&\le&\tau\|z^{k+1}-z^*\|\|\hat{F}(z^k)-\hat{F}(z^{k-1})\|\nonumber\\
&\overset{\eqref{sto-Lip-bd}}{\le}&\tau L\|z^{k+1}-z^*\|\left(\frac{1}{L}(\varepsilon_{z^k}+\varepsilon_{z^{k-1}})+\|z^k-z^{k-1}\|\right)\nonumber\\
&\le& \frac{\tau L}{2}\|z^{k+1}-z^*\|^2+\frac{\tau}{L}(\varepsilon_{z^k}+\varepsilon_{z^{k-1}})^2+\tau L\|z^k-z^{k-1}\|^2\nonumber\\
&\le& \frac{\tau L}{2}\|z^{k+1}-z^*\|^2+\frac{\tau}{L}(\varepsilon_{z^k}+\varepsilon_{z^{k-1}})^2+2\tau L\|z^k-z^{*}\|^2+2\tau L\|z^{k-1}-z^*\|^2. \label{bd-second}
\end{eqnarray}

Combining the bounds in \eqref{bd-last-final} and \eqref{bd-second} with \eqref{bd-first-step} and multiplying both sides by 2, we have
\begin{eqnarray}
&&(1-4|\gamma-\beta|-\tau L)\|z^{k+1}-z^*\|^2\nonumber\\
&\le& \left(1+4\gamma+6|\gamma-\beta|+4\tau L)\right)\|z^k-z^*\|^2\nonumber\\
&&+\left(2|\gamma-\beta|+2\gamma+4\tau L\right)\|z^{k-1}-z^*\|^2\nonumber\\
&&+\left(2\alpha^2L^2+2|\gamma-\beta|+2\gamma-1\right)\|z^{k+0.5}-z^k\|^2\nonumber\\
&&+2\alpha^2(\varepsilon_{z^k}+\varepsilon_{z^{k+0.5}})^2+\frac{2\tau}{L}(\varepsilon_{z^k}+\varepsilon_{z^{k-1}})^2\nonumber\\
&&-2\alpha(z^{k+0.5}-z^*)^\top\hat{F}(z^{k+0.5})+2(\eta-\alpha)(z^{k+1}-z^{k+0.5})^\top\hat{F}(z^k).\nonumber
\end{eqnarray}

Let us now take expectation on both sides. Noting $d_{k+1}=\mathbb{E}\left[\|z^{k+1}-z^*\|^2\right]$, $d_{k}=\mathbb{E}\left[\|z^{k}-z^*\|^2\right]$, and $d_{k-1}=\mathbb{E}\left[\|z^{k-1}-z^*\|^2\right]$, and noting that $\mathbb{E}[\varepsilon_{z}^2]\le\sigma^2$ for all $z\in\mathcal{Z}$ by Assumption \eqref{oracle-var}, we obtain
\begin{eqnarray}
&&(1-4|\gamma-\beta|-\tau L)d_{k+1}\nonumber\\
&\le& \left(1+4\gamma+6|\gamma-\beta|+4\tau L\right)d_k
+\left(2|\gamma-\beta|+2\gamma+4\tau L\right)d_{k-1}\nonumber\\
&&+\left(2\alpha^2L^2+2|\gamma-\beta|+2\gamma-1\right)\mathbb{E}\left[\|z^{k+0.5}-z^k\|^2\right]
+8\left(\alpha^2+\frac{\tau}{L}\right)\sigma^2\nonumber\\
&&-2\alpha\mathbb{E}\left[(z^{k+0.5}-z^*)^\top\hat{F}(z^{k+0.5})\right]+2(\eta-\alpha)\mathbb{E}\left[(z^{k+1}-z^{k+0.5})^\top\hat{F}(z^k)\right].\quad\label{bd-final-0}
\end{eqnarray}

Notice that
\begin{eqnarray}
&&\mathbb{E}\left[(z^{k+0.5}-z^*)^\top\hat{F}(z^{k+0.5})\right]\nonumber\\
&=&\mathbb{E}\left[(z^{k+0.5}-z^*)^\top F(z^{k+0.5})\right]+\mathbb{E}\left[(z^{k+0.5}-z^*)^\top\left(\hat{F}(z^{k+0.5})-F(z^{k+0.5})\right)\right],\nonumber
\end{eqnarray}
where
\begin{eqnarray}
&&\mathbb{E}\left[(z^{k+0.5}-z^*)^\top F(z^{k+0.5})\right]\nonumber\\
&=&\mathbb{E}\left[(z^{k+0.5}-z^*)^\top \left(F(z^{k+0.5})-F(z^*)\right)+(z^{k+0.5}-z^*)^\top F(z^{*})\right]\nonumber\\
&\ge& \mathbb{E}\left[\mu\|z^{k+0.5}-z^*\|^2\right]\nonumber\\
&\ge& \frac{\mu}{2}d_k-\mu\mathbb{E}\left[\|z^{k+0.5}-z^k\|^2\right],\nonumber
\end{eqnarray}
and
\begin{eqnarray}
&&\mathbb{E}\left[(z^{k+0.5}-z^*)^\top\left(\hat{F}(z^{k+0.5})-F(z^{k+0.5})\right)\right]\nonumber\\
&\ge& -\mathbb{E}\left[\|z^{k+0.5}-z^*\|\|\hat{F}(z^{k+0.5})-F(z^{k+0.5})\|\right] \nonumber\\
&=& -\mathbb{E}\left[\mathbb{E}\left[\|z^{k+0.5}-z^*\|\|\hat{F}(z^{k+0.5})-F(z^{k+0.5})\||\xi^{[k]}\right]\right]\nonumber\\
&=& -\mathbb{E}\left[\|z^{k+0.5}-z^*\|\mathbb{E}\left[\|\hat{F}(z^{k+0.5})-F(z^{k+0.5})\||\xi^{[k]}\right]\right] \nonumber\\
&\ge& -\mathbb{E}\left[\|z^{k+0.5}-z^*\|\cdot \delta\right]\nonumber\\
&\ge & -\delta D. \nonumber
\end{eqnarray}
Further note that we have denoted $\xi^{[k]}=(\xi^0,\xi^{0.5},\xi^1,\xi^{1.5},...,\xi^{k-0.5},\xi^k)$ to be the collection of random vectors sampled up until the iterate $z^k$. Therefore, $z^{k+0.5}$ is a known vector given $\xi^{[k]}$.

Putting the above two bounds back into \eqref{bd-final-0}, we arrive at the desired bound:
\begin{eqnarray}
&&(1-4|\gamma-\beta|-\tau L)d_{k+1}\nonumber\\
&\le& \left(1+4\gamma+6|\gamma-\beta|+4\tau L-\alpha\mu\right)d_k\nonumber\\
&&+\left(2|\gamma-\beta|+2\gamma+4\tau L\right)d_{k-1}\nonumber\\
&&+\left(2\alpha^2L^2+2|\gamma-\beta|+2\gamma+2\alpha\mu-1\right)\mathbb{E}\left[\|z^{k+0.5}-z^k\|^2\right]\nonumber\\
&&+8\left(\alpha^2+\frac{\tau}{L}\right)\sigma^2+2\alpha\delta D\nonumber\\
&&+2(\eta-\alpha)\mathbb{E}\left[(z^{k+1}-z^{k+0.5})^\top\hat{F}(z^k)\right].\quad\nonumber
\end{eqnarray}

\label{proof:sto-exp-relation}
\subsection{Proof of Theorem \ref{th:sto-exp-result}}

By condition \eqref{sto-exp-cond-t}, we have $t_3<1$. Let us start with divide both sides of \eqref{sto-exp-relation-2} with $1-t_3$:
\begin{eqnarray}
d_{k+1}&\le& \left(1-\frac{t_1-t_3}{1-t_3}\right)d_k+\frac{t_2}{1-t_3}d_{k-1}+\frac{8\left(\alpha^2+\frac{\tau}{L}\right)\sigma^2}{1-t_3}+\frac{2\alpha\delta D}{1-t_3}\nonumber\\
&=& (1-a)d_k+b\cdot d_{k-1}+c\cdot \sigma^2+d\cdot\delta D.\label{bound-abcd}
\end{eqnarray}
Note that we have $1>a>b$ by condition \eqref{sto-exp-cond-t}. It is elementary to verify that
\[
b\le \left(1-\frac{a-b}{2}\right)\cdot \frac{a+b}{2},
\]
and by rearranging terms in \eqref{bound-abcd}, we have the following
\begin{eqnarray}
d_{k+1}+\frac{a+b}{2}d_k&\le& \left(1-\frac{a-b}{2}\right)d_k+b\cdot d_{k-1}+c\cdot \sigma^2+d\cdot \delta D\nonumber\\
&\le& \left(1-\frac{a-b}{2}\right)\left(d_k+\frac{a+b}{2}d_{k-1}\right)+c\cdot \sigma^2+d\cdot \delta D\nonumber
\end{eqnarray}
A recursive argument yields the following result:
\begin{eqnarray}
&&d_{k+1}+\frac{a+b}{2}d_k\nonumber\\
&\le&\left(1-\frac{a-b}{2}\right)^{k+1}\left(d_0+\frac{a+b}{2}d_{-1}\right)+\left(c\cdot \sigma^2+d\cdot \delta D\right)\cdot\sum\limits_{i=0}^k\left(1-\frac{a-b}{2}\right)^i\nonumber\\
&\le& \left(1-\frac{a-b}{2}\right)^{k+1}\cdot \frac{2+a+b}{2}\|z^0-z^*\|^2+\left(c\cdot \sigma^2+d\cdot \delta D\right)\cdot\frac{2}{a-b}.\nonumber
\end{eqnarray}
Note that $d_0=d_{-1}=\|z^0-z^*\|^2$. The statement in Theorem \ref{th:sto-exp-result} follows by letting $q=\frac{2}{a-b}=\frac{2(1-t_3)}{t_1-t_2-t_3}$.
% and using $\sigma^2=o\left(\left(1-\frac{1}{q}\right)^K\right)$.

\label{proof:sto-exp-result}
\subsection{Proof of Proposition \ref{prop:sto-exp-example}}
\label{proof:sto-exp-example}

For the choice of parameters
\[
(\alpha,\beta,\gamma,\eta,\tau)=\left(\frac{1}{4L},\frac{1}{64\kappa},\frac{1}{64\kappa},\frac{1}{4L},\frac{1}{64L\kappa}\right),
\]
we have
\[
(t_1,t_2,t_3)=\left(\frac{1}{8\kappa},\frac{3}{32\kappa},\frac{1}{64\kappa}\right)
\]
by the relation \eqref{sto-exp-t}. Additionally,
\[
2\alpha^2L^2+2|\gamma-\beta|+2\gamma+2\alpha\mu-1=\frac{1}{8}+\frac{1}{32\kappa}+\frac{1}{2\kappa}-1<0.
\]
Therefore, both conditions \eqref{sto-exp-cond-1} and \eqref{sto-exp-cond-t} are satisfied.

Now, from \eqref{sto-exp-relation}, we have
\begin{eqnarray}
\left(1-\frac{1}{64\kappa}\right)d_{k+1}&\le& \left(1-\frac{1}{8\kappa}\right)d_k+\frac{3}{32\kappa}d_{k-1}+8\left(\frac{1}{16L^2}+\frac{1}{64L^2\kappa}\right)\sigma^2+\frac{\delta D}{2L}\nonumber\\
&\le& \left(1-\frac{1}{8\kappa}\right)d_k+\frac{3}{32\kappa}d_{k-1}+\frac{5\sigma^2}{8L^2}+\frac{\delta D}{2L}.\nonumber
\end{eqnarray}

Divide both sides with $1-\frac{1}{64\kappa}$ and note that $\left(1-\frac{1}{64\kappa}\right)^{-1}\le \frac{64}{63}$, we have:
\begin{eqnarray}
d_{k+1} &\le& \frac{1-\frac{1}{8\kappa}}{1-\frac{1}{64\kappa}}d_k+\frac{2}{21\kappa}d_{k-1}+\frac{40\sigma^2}{63L^2}+\frac{32\delta D}{63L}\nonumber\\
&=& \left(1-\frac{\frac{7}{64\kappa}}{1-\frac{1}{64\kappa}}\right)d_k+\frac{2}{21\kappa}d_{k-1}+\frac{40\sigma^2}{63L^2}+\frac{32\delta D}{63L}\nonumber\\
&\le& \left(1-\frac{7}{64\kappa}\right)d_k+\frac{2}{21\kappa}d_{k-1}+\frac{40\sigma^2}{63L^2}+\frac{32\delta D}{63L}.\nonumber
\end{eqnarray}
We can move a part of $d_k$ to the LHS and form the following:
\begin{eqnarray}
&&d_{k+1}+\frac{27}{256\kappa}d_k\nonumber\\
&\le& \left(1-\frac{1}{256\kappa}\right)d_k+\frac{2}{21\kappa}d_{k-1}+\frac{40\sigma^2}{63L^2}+\frac{32\delta D}{63L}\nonumber\\
&\le& \left(1-\frac{1}{256\kappa}\right)\left(d_k+\frac{27}{256}d_{k-1}\right)+\frac{40\sigma^2}{63L^2}+\frac{32\delta D}{63L}\nonumber\\
&\le& \left(1-\frac{1}{256\kappa}\right)^{k+1}\left(d_0+\frac{27}{256}d_{-1}\right)+\left(\frac{40\sigma^2}{63L^2}+\frac{32\delta D}{63L}\right)\cdot\sum\limits_{i=0}^{k}\left(1-\frac{1}{256\kappa}\right)^i\nonumber\\
&=& \left(1-\frac{1}{256\kappa}\right)^{k+1}\cdot\frac{283}{256}\|z^0-z^*\|^2+\left(\frac{40\sigma^2}{63L^2}+\frac{32\delta D}{63L}\right)\cdot\frac{1-\left(1-\frac{1}{256\kappa}\right)^{k+1}}{\frac{1}{256\kappa}}\nonumber\\
&\le& \left(1-\frac{1}{256\kappa}\right)^{k+1}\cdot\frac{283}{256}\|z^0-z^*\|^2+\left(\frac{40\sigma^2}{63L^2}+\frac{32\delta D}{63L}\right)\cdot256\kappa.\nonumber
\end{eqnarray}
Note that $d_0=d_{-1}=\|z^0-z^*\|^2$. Finally, the LHS of the above inequality can be lower bounded by $d_{k+1}$, thus completing the proof.

\subsection{Proof of Lemma \ref{lem:sto-exm-relation}}
\label{proof:sto-exm-relation}
We start by using the 1-co-coerciveness of the projection operator $P_{\mathcal{Z}}(\cdot)$:
\begin{eqnarray}
&&\|z^{k+1}-z^*\|^2\nonumber\\
&=& \|P_{\mathcal{Z}}\left(z^k-\alpha \hat{F}(z^k)+\gamma(z^k-z^{k-1})-\tau\left(\hat{F}(z^k)-\hat{F}(z^{k-1})\right)\right)-P_{\mathcal{Z}}\left(z^*-\alpha F(z^*)\right)\|^2\nonumber\\
&\le& (z^{k+1}-z^*)^\top\left(z^k-\alpha \hat{F}(z^k)+\gamma(z^k-z^{k-1})-\tau\left(\hat{F}(z^k)-\hat{F}(z^{k-1})\right)-\left(z^*-\alpha F(z^*)\right)\right)\nonumber\\
&=& (z^{k+1}-z^*)^\top\left((z^k-z^*)-\alpha \left(\hat{F}(z^k)-F(z^*)\right)+\gamma(z^k-z^{k-1})-\tau\left(\hat{F}(z^k)-\hat{F}(z^{k-1})\right)\right). \nonumber\\
&&\label{sto-OGDA-hy-first-step}
\end{eqnarray}

Next, let us bound the above four terms separately:
\begin{eqnarray}
(z^{k+1}-z^*)^\top(z^k-z^*)=\frac{1}{2}\left(\|z^{k+1}-z^*\|^2+\|z^k-z^*\|^2-\|z^{k+1}-z^k\|^2\right),\label{first-step-1}
\end{eqnarray}
and
\begin{eqnarray}
(z^{k+1}-z^*)^\top\left(\hat{F}(z^k)-F(z^*)\right)&=& (z^{k+1}-z^*)^\top\left(\hat{F}(z^k)-\hat{F}(z^{k+1})+\hat{F}(z^{k+1})-F(z^*)\right) \label{first-step-2-1}
\end{eqnarray}
where
\begin{eqnarray}
&&(z^{k+1}-z^*)^\top\left(\hat{F}(z^{k+1})-F(z^*)\right)\nonumber\\
&=&(z^{k+1}-z^*)^\top\left(\hat{F}(z^{k+1})-F(z^{k+1})+F(z^{k+1})-F(z^*)\right)\nonumber\\
&\ge& (z^{k+1}-z^*)^\top\left(\hat{F}(z^{k+1})-F(z^{k+1})\right)+\mu\|z^{k+1}-z^*\|^2,\label{first-step-2-2}
\end{eqnarray}
and
\begin{eqnarray}
(z^{k+1}-z^*)^\top(z^k-z^{k-1})&\le& \frac{\gamma}{2}\left(\|z^{k+1}-z^*\|^2+\|z^k-z^{k-1}\|^2\right),\label{first-step-3}
\end{eqnarray}
and
\begin{eqnarray}
&&-\tau(z^{k+1}-z^*)^\top\left(\hat{F}(z^k)-\hat{F}(z^{k-1})\right)\nonumber\\
&=& -\tau(z^{k+1}-z^k+z^k-z^*)^\top\left(\hat{F}(z^k)-\hat{F}(z^{k-1})\right)\nonumber\\
&\le& \tau\|z^{k+1}-z^k\|\|\hat{F}(z^k)-\hat{F}(z^{k-1})\|-\tau(z^k-z^*)^\top\left(\hat{F}(z^k)-\hat{F}(z^{k-1})\right)\nonumber\\
&\le& \frac{1}{4}\|z^{k+1}-z^k\|^2+\tau^2\|\hat{F}(z^k)-\hat{F}(z^{k-1})\|^2-\tau(z^k-z^*)^\top\left(\hat{F}(z^k)-\hat{F}(z^{k-1})\right),\label{first-step-4-1}
\end{eqnarray}
where
\begin{eqnarray}
\tau^2\|\hat{F}(z^k)-\hat{F}(z^{k-1})\|^2&\overset{\eqref{sto-Lip-bd}}{\le} & 2\tau^2(\varepsilon_{z^k}+\varepsilon_{z^{k-1}})^2+2\tau^2L^2\|z^k-z^{k-1}\|^2.\label{first-step-4-2}
\end{eqnarray}

Putting the bounds \eqref{first-step-1}-\eqref{first-step-4-2} back to \eqref{sto-OGDA-hy-first-step}, we get:
\begin{eqnarray}
&&\left(\frac{1}{2}+\alpha\mu-\frac{\gamma}{2}\right)\|z^{k+1}-z^*\|^2+\alpha(z^{k+1}-z^*)^\top\left(\hat{F}(z^k)-\hat{F}(z^{k+1})\right)+\frac{1}{4}\|z^{k+1}-z^k\|^2\nonumber\\
&\le& \frac{1}{2}\|z^k-z^*\|^2+\tau(z^k-z^*)^\top\left(\hat{F}(z^{k-1})-\hat{F}(z^{k})\right)+\left(2\tau^2L^2+\frac{\gamma}{2}\right)\|z^k-z^{k-1}\|^2\nonumber\\
&&+2\tau^2(\varepsilon_{z^k}+\varepsilon_{z^{k-1}})^2-\alpha(z^{k+1}-z^*)^\top\left(\hat{F}(z^{k+1})-F(z^{k+1})\right).\nonumber
\end{eqnarray}

Taking expectation on both sides gives us
\begin{eqnarray}
&&\mathbb{E}\left[\left(\frac{1}{2}+\alpha\mu-\frac{\gamma}{2}\right)\|z^{k+1}-z^*\|^2+\alpha(z^{k+1}-z^*)^\top\left(\hat{F}(z^k)-\hat{F}(z^{k+1})\right)+\frac{1}{4}\|z^{k+1}-z^k\|^2\right]\nonumber\\
&\le& \mathbb{E}\left[\frac{1}{2}\|z^k-z^*\|^2+\tau(z^k-z^*)^\top\left(\hat{F}(z^{k-1})-\hat{F}(z^{k})\right)+\left(2\tau^2L^2+\frac{\gamma}{2}\right)\|z^k-z^{k-1}\|^2\right]\nonumber\\
&&+8\tau^2\sigma^2-\mathbb{E}\left[\alpha(z^{k+1}-z^*)^\top\left(\hat{F}(z^{k+1})-F(z^{k+1})\right)\right].\label{sto-OGDA-hy-second-step}
\end{eqnarray}
Note that
\begin{eqnarray}
&&-\mathbb{E}\left[\alpha(z^{k+1}-z^*)^\top\left(\hat{F}(z^{k+1})-F(z^{k+1})\right)\right]\nonumber\\
&\le&\alpha\mathbb{E}\left[\|z^{k+1}-z^*\|\|\hat{F}(z^{k+1})-F(z^{k+1})\|\right]\nonumber\\
&=&\alpha\mathbb{E}\left[\mathbb{E}\left[\left(\|z^{k+1}-z^*\|\|\hat{F}(z^{k+1})-F(z^{k+1})\|\right)|\xi^{[k]}\right]\right]\nonumber\\
&=&\alpha\mathbb{E}\left[\|z^{k+1}-z^*\|\cdot\mathbb{E}\left[\left(\|\hat{F}(z^{k+1})-F(z^{k+1})\|\right)|\xi^{[k]}\right]\right]\nonumber\\
&\le& \alpha\mathbb{E}\left[\delta\|z^{k+1}-z^*\|\right]\nonumber\\
&\le& \alpha\mathbb{E}\left[\frac{\delta^2}{2\mu}+\frac{\mu}{2}\|z^{k+1}-z^*\|^2\right].\nonumber
\end{eqnarray}
Here we define $\xi^{[k]}=(\xi^0,\xi^1,...,\xi^k)$ to be the collection of random vectors sampled up until the iterate $z^k$, and $z^{k+1}$ is known given $\xi^{[k]}$.

Therefore, \eqref{sto-OGDA-hy-second-step} becomes
\begin{eqnarray}
&&\mathbb{E}\left[\left(\frac{1}{2}+\frac{\alpha\mu}{2}-\frac{\gamma}{2}\right)\|z^{k+1}-z^*\|^2+\alpha(z^{k+1}-z^*)^\top\left(\hat{F}(z^k)-\hat{F}(z^{k+1})\right)+\frac{1}{4}\|z^{k+1}-z^k\|^2\right]\nonumber\\
&\le& \mathbb{E}\left[\frac{1}{2}\|z^k-z^*\|^2+\tau(z^k-z^*)^\top\left(\hat{F}(z^{k-1})-\hat{F}(z^{k})\right)+\left(2\tau^2L^2+\frac{\gamma}{2}\right)\|z^k-z^{k-1}\|^2\right]\nonumber\\
&&+8\tau^2\sigma^2+\frac{\alpha\delta^2}{2\mu},\nonumber
\end{eqnarray}
completing the proof.

\subsection{Proof of Theorem \ref{th:sto-exm-result}}
\label{proof:sto-exm-result}

Continuing from \eqref{sto-exm-V-relation}, we have:
\begin{eqnarray*}
V_k&\le& \left(1+\frac{\theta}{\kappa}\right)^{-k}V_0+\sum\limits_{i=1}^{k}\left(1+\frac{\theta}{\kappa}\right)^{-i}\cdot\left(8\tau^2\sigma^2+\frac{\alpha\delta^2}{2\mu}\right)\nonumber\\
&=& \frac{1}{2}\left(1+\frac{\theta}{\kappa}\right)^{-k}\|z^0-z^*\|^2+\frac{1-\left(1+\frac{\theta}{\kappa}\right)^{-k}}{\frac{\theta}{\kappa}}\cdot\left(8\tau^2\sigma^2+\frac{\alpha\delta^2}{2\mu}\right),\nonumber
\end{eqnarray*}
where we use $z^{-1}=z^0$ for $V_0$.

Finally, with the following bound:
\begin{eqnarray*}
&&\tau(z^k-z^*)^\top\left(\hat{F}(z^{k-1})-\hat{F}(z^{k})\right)\\
&\ge& -\tau \|z^k-z^*\|\|\hat{F}(z^{k-1})-\hat{F}(z^{k})\|\\
&\ge& -\frac{1}{4}\|z^k-z^*\|^2-\tau^2\|\hat{F}(z^{k-1})-\hat{F}(z^{k})\|^2\\
&\overset{\eqref{sto-Lip-bd}}{\ge}& -\frac{1}{4}\|z^k-z^*\|^2-\tau^2\left(2L^2\|z^{k-1}-z^k\|^2+2(\varepsilon_{z^{k-1}}+\varepsilon_{z^k})^2\right),
\end{eqnarray*}
we can lower bound $V_k$ as
\[
\frac{1}{4}\mathbb{E}\left[\|z^k-z^*\|^2\right]-8\tau^2\sigma^2\le V_k.
\]
Therefore
\begin{eqnarray}
&&\mathbb{E}\left[\|z^k-z^*\|^2\right]\nonumber\\
&\le& 2\left(1+\frac{\theta}{\kappa}\right)^{-k}\|z^0-z^*\|^2+\frac{4\kappa}{\theta}\cdot\left(1-\left(1+\frac{\theta}{\kappa}\right)^{-k}\right)\cdot\left(8\tau^2\sigma^2+\frac{\alpha\delta^2}{2\mu}\right)+32\tau^2\sigma^2\nonumber\\
&\le& 2\left(1+\frac{\theta}{\kappa}\right)^{-k}\|z^0-z^*\|^2+\left(\frac{\kappa}{\theta}+1\right)\cdot 32\tau^2\sigma^2+\frac{2\kappa\alpha\delta^2}{\theta\mu}.\label{sto-exm-final-ineq}
\end{eqnarray}
The statement in Theorem \ref{th:sto-exm-result} follows by noting $\left(1+\frac{\theta}{\kappa}\right)^{-1}=1-\frac{\theta}{\kappa+\theta}$.
% and using the condition $\sigma^2=o\left(\left(1-\frac{1}{\kappa}\right)^K\right)$.

\subsection{Proof of Lemma \ref{lem:sto-0-grad}}
\label{proof:sto-0-grad}
We will derive the first bound in \eqref{SZG-var}; the second bound is similar and will be omitted. %while omitting the similar proof for the second bound.

Notice that
\begin{eqnarray*}
&&\mathbb{E}_{\xi,u}\left[\|F_{\rho_x}(x,y,\xi,u)\|^2\right]\\
&=& \mathbb{E}_{\xi}\left[\mathbb{E}_{u}\left[\|F_{\rho_x}(x,y,\xi,u)\|^2\right]\right]\\
&\overset{\eqref{smooth-x-sm}}{\le}& \mathbb{E}_{\xi}\left[2n\|\nabla_x\hat{f}(x,y,\xi)\|^2\right]+\frac{\rho_x^2L^2n^2}{2}\\
&=& 2n\left(\mathbb{E}_{\xi}\left[\|\nabla_xf(x,y)\|^2+2\nabla_xf(x,y)^\top\left(\nabla_x\hat{f}(x,y,\xi)-\nabla_xf(x,y)\right)+\|\nabla_x\hat{f}(x,y,\xi)-\nabla_xf(x,y)\|^2\right]\right)\\
&&+\frac{\rho_x^2L^2n^2}{2}\\
&\overset{\eqref{so-x-mean}}{=}& 2n\left(\|\nabla_xf(x,y)\|^2+\mathbb{E}_{\xi}\left[\|\nabla_x\hat{f}(x,y,\xi)-\nabla_xf(x,y)\|^2\right]\right)+\frac{\rho_x^2L^2n^2}{2}\\
&\le& 2n\left(M^2+\sigma^2\right)+\frac{\rho_x^2L^2n^2}{2}.
\end{eqnarray*}
Further note that
\begin{eqnarray*}
&& \mathbb{E}_{\xi,u}\left[\|F_{\rho_x}(x,y,\xi,u)-\nabla_xf_{\rho_{x}}(x,y)\|^2\right]\\
&=&\mathbb{E}_{\xi,u}\left[\|F_{\rho_x}(x,y,\xi,u)\|^2-2F_{\rho_x}(x,y,\xi,u)^\top\nabla_xf_{\rho_x}(x,y)+\|\nabla_xf_{\rho_{x}}(x,y)\|^2\right]\\
&\overset{\eqref{SZG-mean}}{=}&\mathbb{E}_{\xi,u}\left[\|F_{\rho_x}(x,y,\xi,u)\|^2-\|\nabla_xf_{\rho_{x}}(x,y)\|^2\right] \\
&\le& \mathbb{E}_{\xi,u}\left[\|F_{\rho_x}(x,y,\xi,u)\|^2\right]\\
&\le& 2n(M^2+\sigma^2)+\frac{\rho_x^2L^2n^2}{2},
\end{eqnarray*}
completing the proof for \eqref{SZG-var}.

\subsection{Proof of Lemma \ref{lem:sto-0-exp-relation}}
\label{proof:sto-0-exp-relation}

The logic line of the proof for this lemma is very similar to the proof in Appendix \ref{proof:sto-exp-relation}, with the stochastic mapping $\hat{F}(z^k)$ replaced by the stochastic zeroth-order gradient $\hat{F}^k_{\rho}(z^k)$. Therefore, we shall refrain from repeating similar analysis, but highlight the main differences instead.
First, for $F(z)=\begin{pmatrix}\nabla_xf(x,y)\\-\nabla_yf(x,y)\end{pmatrix}$, we shall have
\begin{eqnarray}
\|F(z)-F(z')\|\le L\|z-z'\|,\quad\forall z,z'\in\mathcal{Z}\label{lip-minmax}
\end{eqnarray}
where $L=2\cdot\max(L_x,L_y,L_{xy})$, because %from the following inequalities
\begin{eqnarray}
&&\|F(z)-F(z')\|^2\nonumber\\
&=& \|\nabla_xf(x,y)-\nabla_xf(x',y')\|^2+\|\nabla_yf(x,y)-\nabla_yf(x',y')\|^2\nonumber\\
&\le& 2L_x^2\|x-x'\|^2+2L_{xy}^2\|y-y'\|^2+2L_y^2\|y-y'\|^2+2L_{xy}^2\|x-x'\|^2\nonumber\\
&\le& L^2(\|x-x'\|^2+\|y-y'\|^2)=L^2\|z-z'\|^2.\nonumber
\end{eqnarray}

Next, by denoting $\varepsilon_{z^k}=\|\hat{F}^k_{\rho}(z^k)-F_{\rho}(z^k)\|$, we have
\begin{eqnarray}
&&\|\hat{F}^k_{\rho}(z^k)-\hat{F}^{k+0.5}_{\rho}(z^{k+0.5})\|\nonumber\\
&=& \|\hat{F}^k_{\rho}(z^k)-F_{\rho}(z^k)+F_\rho(z^{k+0.5})-\hat{F}_{\rho}(z^{k+0.5})+F_{\rho}(z^k)-F_{\rho}(z^{k+0.5})\|\nonumber\\
&\le& \varepsilon_{z^k}+\varepsilon_{z^{k+0.5}}+\|F_{\rho}(z^k)-F_{\rho}(z^{k+0.5})\|\nonumber\\
&=&\varepsilon_{z^k}+\varepsilon_{z^{k+0.5}}+\|F_{\rho}(z^{k})-F(z^k)+F(z^{k+0.5})-F_{\rho}(z^{k+0.5})+F(z^k)-F(z^{k+0.5})\| \nonumber\\
&\overset{\eqref{smooth-grad-bd},\eqref{lip-minmax}}{\le}&  \varepsilon_{z^k}+\varepsilon_{z^{k+0.5}}+L\sqrt{\rho_x^2n^2+\rho_y^2m^2}+L\|z^k-z^{k+0.5}\| . \label{Lip-szo-grad}
\end{eqnarray}

Therefore, for a similar bound as in \eqref{sto-exp-lip-bd-1}, we have
\begin{eqnarray}
&&\alpha\|z^{k+1}-z^{k+0.5}\|\|\hat{F}^k_{\rho}(z^k)-\hat{F}^{k+0.5}_{\rho}(z^{k+0.5})\|\nonumber\\
&\le& \frac{1}{2}\left(\|z^{k+1}-z^{k+0.5}\|^2+\alpha^2\|\hat{F}^k_{\rho}(z^k)-\hat{F}^{k+0.5}_{\rho}(z^{k+0.5})\|^2\right)\nonumber\\
&\le& \frac{1}{2}\|z^{k+1}-z^{k+0.5}\|^2+2\alpha^2(\varepsilon_{z^k}+\varepsilon_{z^{k+0.5}})^2+2\alpha^2L^2(\rho_x^2n^2+\rho_y^2m^2)+\alpha^2L^2\|z^k-z^{k+0.5}\|^2.\nonumber
\end{eqnarray}

For another similar bound as in \eqref{bd-second}, we have
\begin{eqnarray}
&&-\tau(z^{k+1}-z^*)^\top\left(\hat{F}^k_{\rho}(z^k)-\hat{F}^{k-1}_{\rho}(z^{k-1})\right)\nonumber\\
&\le&\tau\|z^{k+1}-z^*\|\|\hat{F}^k_{\rho}(z^k)-\hat{F}^{k-1}_{\rho}(z^{k-1})\|\nonumber\\
&\le&\tau L\|z^{k+1}-z^*\|\left(\frac{1}{L}(\varepsilon_{z^k}+\varepsilon_{z^{k-1}})+\sqrt{\rho_x^2n^2+\rho_y^2m^2}+\|z^k-z^{k-1}\|\right)\nonumber\\
&\le& \frac{\tau L}{2}\|z^{k+1}-z^*\|^2+\frac{2\tau}{L}(\varepsilon_{z^k}+\varepsilon_{z^{k-1}})^2+2\tau L(\rho_x^2n^2+\rho_y^2m^2)+\tau L\|z^k-z^{k-1}\|^2\nonumber\\
&\le& \frac{\tau L}{2}\|z^{k+1}-z^*\|^2+\frac{2\tau}{L}(\varepsilon_{z^k}+\varepsilon_{z^{k-1}})^2+2\tau L(\rho_x^2n^2+\rho_y^2m^2)\nonumber\\
&&+2\tau L\|z^k-z^{*}\|^2+2\tau L\|z^{k-1}-z^*\|^2.\nonumber
\end{eqnarray}

Therefore, we reach the bound that
\begin{eqnarray}
&&(1-4|\gamma-\beta|-\tau L)\|z^{k+1}-z^*\|^2\nonumber\\
&\le& \left(1+4\gamma+6|\gamma-\beta|+4\tau L)\right)\|z^k-z^*\|^2\nonumber\\
&&+\left(2|\gamma-\beta|+2\gamma+4\tau L\right)\|z^{k-1}-z^*\|^2\nonumber\\
&&+\left(2\alpha^2L^2+2|\gamma-\beta|+2\gamma-1\right)\|z^{k+0.5}-z^k\|^2\nonumber\\
&&+4\alpha^2(\varepsilon_{z^k}+\varepsilon_{z^{k+0.5}})^2+\frac{4\tau}{L}(\varepsilon_{z^k}+\varepsilon_{z^{k-1}})^2\nonumber\\
&&+4(\alpha^2L^2+\tau L)(\rho_x^2n^2+\rho_y^2m^2)\nonumber\\
&&-2\alpha(z^{k+0.5}-z^*)^\top\hat{F}^{k+0.5}_{\rho}(z^{k+0.5})+2(\eta-\alpha)(z^{k+1}-z^{k+0.5})^\top\hat{F}^k_{\rho}(z^k).\nonumber
\end{eqnarray}

By \eqref{batch-var}, we have $\mathbb{E}[\varepsilon^2_{z^k}]\le\frac{2\tilde{\sigma}^2}{t_k}$. Taking expectation on both sides for the above inequality, we have
\begin{eqnarray}
&&(1-4|\gamma-\beta|-\tau L)d_{k+1}\nonumber\\
&\le& \left(1+4\gamma+6|\gamma-\beta|+4\tau L)\right)d_k\nonumber\\
&&+\left(2|\gamma-\beta|+2\gamma+4\tau L\right)d_{k-1}\nonumber\\
&&+\left(2\alpha^2L^2+2|\gamma-\beta|+2\gamma-1\right)\mathbb{E}\left[\|z^{k+0.5}-z^k\|^2\right]\nonumber\\
&&+16\tilde{\sigma}^2\left(\left(\alpha^2+\frac{\tau}{L}\right)\frac{1}{t_k}+\frac{\alpha^2}{t_{k+0.5}}+\frac{\tau}{Lt_{k-1}}\right)+4(\alpha^2L^2+\tau L)(\rho_x^2n^2+\rho_y^2m^2)\nonumber\\
&&+\mathbb{E}\left[-2\alpha(z^{k+0.5}-z^*)^\top\hat{F}^{k+0.5}_{\rho}(z^{k+0.5})+2(\eta-\alpha)(z^{k+1}-z^{k+0.5})^\top\hat{F}^k_{\rho}(z^k)\right].\label{sto-0-exp-bd-1}
\end{eqnarray}

Note that
\begin{eqnarray}
&& \mathbb{E}\left[(z^{k+0.5}-z^*)^\top\hat{F}^{k+0.5}_{\rho}(z^{k+0.5})\right] \nonumber\\
&=& \mathbb{E}\left[(z^{k+0.5}-z^*)^\top F(z^{k+0.5})\right]+\mathbb{E}\left[(z^{k+0.5}-z^*)^\top\left(\hat{F}^{k+0.5}_{\rho}(z^{k+0.5})-F(z^{k+0.5})\right)\right],\nonumber
\end{eqnarray}
where
\[
\mathbb{E}\left[(z^{k+0.5}-z^*)^\top F(z^{k+0.5})\right]\ge \mathbb{E}\left[\mu\|z^{k+0.5}-z^*\|^2\right]\ge \frac{\mu d_k}{2}-\mu\mathbb{E}\left[\|z^{k+0.5}-z^k\|^2\right].
\]
Let us denote
\begin{eqnarray}
&\xi^k_{[t_k]} := (\xi^k_1,...,\xi^k_{t_k}),& w^k_{[t_k]} := (w^k_1,...,w^k_{t_k})\nonumber\\
&\xi^{[k]} := (\xi^1_{[t_1]},...,\xi^k_{[t_k]}),& w^{[k]} := (w^1_{[t_1]},...,w^k_{[t_k]})\nonumber
\end{eqnarray}
as the collection of all random vectors at iteration $k$ and the collection of all such random vectors from iteration $1$ to $k$ respectively. Notice that for the given $(\xi^{[k]},w^{[k]})$, $z^{k+0.5}$ is a deterministic vector, we then have the following bound
\begin{eqnarray*}
&& \mathbb{E}\left[(z^{k+0.5}-z^*)^\top\left(\hat{F}^{k+0.5}_{\rho}(z^{k+0.5})-F(z^{k+0.5})\right)\right]\\
&=& \mathbb{E}\left[\mathbb{E}\left[(z^{k+0.5}-z^*)^\top\left(\hat{F}^{k+0.5}_{\rho}(z^{k+0.5})-F(z^{k+0.5})\right)\left|\xi^{[k]},w^{[k]}\right.\right]\right]\\
&\overset{\eqref{batch-mean}}{=}& \mathbb{E}\left[(z^{k+0.5}-z^*)^\top\left(F_{\rho}(z^{k+0.5})-F(z^{k+0.5})\right)\right]\\
&\ge& -\mathbb{E}\left[\|z^{k+0.5}-z^*\|\left\|F_{\rho}(z^{k+0.5})-F(z^{k+0.5})\right\|\right]\\
&\overset{\eqref{smooth-grad-bd}}{\ge}& -\frac{L\sqrt{\rho_x^2n^2+\rho_y^2m^2}}{2}\mathbb{E}\left[\|z^{k+0.5}-z^*\|\right]\\
&\ge& -\frac{LD\sqrt{\rho_x^2n^2+\rho_y^2m^2}}{2},
\end{eqnarray*}
where in the last inequality we utilize the boundedness assumption of $\mathcal{Z}=\mathcal{X}\times\mathcal{Y}$ and denote $D=\max\limits_{z,z'\in\mathcal{Z}}\|z-z'\|$. Combining the above bounds into \eqref{sto-0-exp-bd-1}, the desired result follows.

\subsection{Proof of Proposition \ref{prop:szo-exp-sample}}
\label{proof:prop-szo-exp-sample}

Note the variance is upper bounded by
\[
16\tilde{\sigma}^2\left(\left(\alpha^2+\frac{\tau}{L}\right)\frac{1}{t_k}+\frac{\alpha^2}{t_{k+0.5}}+\frac{\tau}{Lt_{k-1}}\right).
\]
Since we take $t_k=t_{k+0.5}$ and $\frac{1}{t_{k-1}}=\frac{1}{t_k}\left(1-\frac{1}{256\kappa}\right)^{-1}\le\frac{2}{t_k}$, the above upper bound can be written as
\[
\left(\alpha^2+\frac{\tau}{L}\right)\frac{48\tilde{\sigma}^2}{t_k}.
\]

By substituting the specific parameter choice into \eqref{szo-exp-relation}, starting from the last iteration $K$, we have
\begin{eqnarray}
&&\left(1-\frac{1}{64\kappa}\right)d_{K}\nonumber\\
&\le& \left(1-\frac{1}{8\kappa}\right)d_{K-1}+\frac{3}{32\kappa}d_{K-2}+\frac{\tilde{\sigma}^2}{t_{K-1}}\cdot\left(\frac{3}{L^2}+\frac{3}{4L^2\kappa}\right)\nonumber\\
&&+\left(\frac{1}{16\kappa}+\frac{1}{4}\right)(\rho_x^2n^2+\rho_y^2m^2)+\frac{D\sqrt{\rho_x^2n^2+\rho_y^2m^2}}{4}\nonumber\\
&\le& \left(1-\frac{1}{8\kappa}\right)d_{K-1}+\frac{3}{32\kappa}d_{K-2}+\frac{15\tilde{\sigma}^2}{4t_{K-1}L^2}+\frac{5(\rho_x^2n^2+\rho_y^2m^2)}{16}+\frac{D\sqrt{\rho_x^2n^2+\rho_y^2m^2}}{4}.\nonumber\\
&\le& \left(1-\frac{1}{8\kappa}\right)d_{K-1}+\frac{3}{32\kappa}d_{K-2}+\frac{15\tilde{\sigma}^2}{4t_{K-1}L^2}+\left(\frac{5}{16}+\frac{D}{4}\right)\frac{C_1^K}{\kappa}.\nonumber
\end{eqnarray}
In the last inequality, we denote $C_1=1-\frac{1}{256\kappa}$ and $\rho_x=\frac{C_1^K}{\sqrt{2}n\kappa}$, $\rho_y=\frac{C_1^K}{\sqrt{2}m\kappa}$, and use the fact that $C_1^2\le C_1$.

Dividing both sides by $1-\frac{1}{64\kappa}$ and noting that $\left(1-\frac{1}{64\kappa}\right)^{-1}\le \frac{64}{63}$, we obtain
\begin{eqnarray}
d_{K} &\le& \frac{1-\frac{1}{8\kappa}}{1-\frac{1}{64\kappa}}d_{K-1}+\frac{2}{21\kappa}d_{K-2}+\frac{80\tilde{\sigma}^2}{21t_{K-1}L^2}+\frac{(20+16D)C_1^K}{63\kappa}\nonumber\\
&=& \left(1-\frac{\frac{7}{64\kappa}}{1-\frac{1}{64\kappa}}\right)d_{K-1}+\frac{2}{21\kappa}d_{K-2}+\frac{80\tilde{\sigma}^2}{21t_{K-1}L^2}+\frac{(20+16D)C_1^K}{63\kappa}\nonumber\\
&\le& \left(1-\frac{7}{64\kappa}\right)d_{K-1}+\frac{2}{21\kappa}d_{K-2}+\frac{80\tilde{\sigma}^2}{21t_{K-1}L^2}+\frac{(20+16D)C_1^K}{63\kappa}.\nonumber\\
&=& \left(1-\frac{7}{64\kappa}\right)d_{K-1}+\frac{2}{21\kappa}d_{K-2}+\frac{c_2}{t_{K-1}}+c_3\cdot \frac{C_1^K}{\kappa}.\nonumber
\end{eqnarray}

By moving a part of $d_{K-1}$ to the LHS, we have:
\begin{eqnarray}
&&d_{K}+\frac{27}{256\kappa}d_{K-1}\nonumber\\
&\le& C_1\cdot d_{K-1}+\frac{2}{21\kappa}d_{K-2}+\frac{c_2}{t_{K-1}}+c_3\cdot \frac{C_1^K}{\kappa}\nonumber\\
&\le& C_1\left(d_{K-1}+\frac{27}{256}d_{K-2}\right)+\frac{c_2}{t_{K-1}}+c_3\cdot \frac{C_1^K}{\kappa}\nonumber\\
&\le& C_1^{K}\left(d_0+\frac{27}{256}d_{-1}\right)+c_2\cdot\sum\limits_{i=1}^{K}t_{K-i}^{-1}C_1^{i-1}+c_3\frac{C_1^K}{\kappa}\cdot\sum\limits_{i=1}^{K}C_1^{i-1}\nonumber\\
&=& C_1^{K}\cdot\frac{283}{256}\|z^0-z^*\|^2+c_2\cdot\sum\limits_{i=1}^{K}t_{K-i}^{-1}C_1^{i-1}+c_3\frac{C_1^K}{\kappa}\cdot\frac{1-C_1^{K}}{\frac{1}{256\kappa}}\nonumber\\
&\le& C_1^{K}\cdot\frac{283}{256}\|z^0-z^*\|^2+c_2\cdot\sum\limits_{i=1}^{K}t_{K-i}^{-1}C_1^{i-1}+256c_3C_1^K.\nonumber
\end{eqnarray}
With the sample size $t_k=K\left(1-\frac{1}{256\kappa}\right)^{-(k+1)}=K\cdot C_1^{-(k+1)}$, we have:
\[
d_K\le \left(1-\frac{1}{256\kappa}\right)^{K}\left(\frac{283}{256}\|z^0-z^*\|^2+c_2+256c_3\right).
\]
It is then straightforward to see that for $K=\mathcal{O}\left(\kappa\ln\left(\frac{1}{\epsilon}\right)\right)$ we have $d_K\le \epsilon$, with the sample complexity given by
\begin{eqnarray}
\sum\limits_{k=0}^{K-1}(t_k+t_{k+0.5})=2\sum\limits_{k=0}^{K-1}t_k=2K\cdot \frac{1-C_1^{-K}}{C_1(1-C_1^{-1})} = \frac{2K(C_1^{-K}-1)}{1-C_1}.\nonumber
\end{eqnarray}

By noticing that a more precise expression of $K=\ln_{C_1^{-1}}\left(\frac{1}{\epsilon}\right)$, we have $C_1^{-K}=\mathcal{O}\left(\frac{1}{\epsilon}\right)$, $1-C_1=\mathcal{O}\left(\frac{1}{\kappa}\right)$, the combined sample complexity is then given by $\mathcal{O}\left(\frac{\kappa^2}{\epsilon}\ln\left(\frac{1}{\epsilon}\right)\right)$.

\subsection{Proof of Lemma \ref{lem:szo-exm-relation}}
\label{proof:szo-exm-relation}

With the similar logic to the proof in Appendix \ref{proof:sto-exm-relation}, we shall focus on the main differences between the two proofs.

Firstly, with the similar derivation to \eqref{Lip-szo-grad}, we have the following bound:
\begin{eqnarray}
&&\tau^2\|\hat{F}^k_{\rho}(z^k)-\hat{F}^{k-1}_{\rho}(z^{k-1})\|^2\nonumber\\
&\le & \tau^2\left(\varepsilon_{z^k}+\varepsilon_{z^{k-1}}+L\sqrt{\rho_x^2n^2+\rho_y^2m^2}+L\|z^k-z^{k-1}\|\right)^2\nonumber\\
&\le& \tau^2\left(4(\varepsilon_{z^k}+\varepsilon_{z^{k-1}})^2+4L^2(\rho_x^2n^2+\rho_y^2m^2)+2L^2\|z^k-z^{k-1}\|^2\right).\nonumber
\end{eqnarray}
By using the variance bound in \eqref{batch-var}, we will reach the next inequality that is similar to the step in \eqref{sto-OGDA-hy-second-step}:
\begin{eqnarray}
&&\mathbb{E}\left[\left(\frac{1}{2}+\alpha\mu-\frac{\gamma}{2}\right)\|z^{k+1}-z^*\|^2+\alpha(z^{k+1}-z^*)^\top\left(\hat{F}^k_{\rho}(z^k)-\hat{F}^{k+1}_{\rho}(z^{k+1})\right)+\frac{1}{4}\|z^{k+1}-z^k\|^2\right]\nonumber\\
&\le& \mathbb{E}\left[\frac{1}{2}\|z^k-z^*\|^2+\tau(z^k-z^*)^\top\left(\hat{F}^{k-1}_{\rho}(z^{k-1})-\hat{F}^{k}_{\rho}(z^{k})\right)+\left(2\tau^2L^2+\frac{\gamma}{2}\right)\|z^k-z^{k-1}\|^2\right]\nonumber\\
&&+16\tau^2\sigma^2\left(\frac{1}{t_k}+\frac{1}{t_{k-1}}\right)+4\tau^2L^2(\rho_x^2n^2+\rho_y^2m^2)\nonumber\\
&&-\mathbb{E}\left[\alpha(z^{k+1}-z^*)^\top\left(\hat{F}^{k+1}_{\rho}(z^{k+1})-F(z^{k+1})\right)\right].\label{szo-exm-bd-1}
\end{eqnarray}

Denote by
\begin{eqnarray}
&\xi^k_{[t_k]}=(\xi^k_1,...,\xi^k_{t_k}),& w^k_{[t_k]}=(w^k_1,...,w^k_{t_k})\nonumber\\
&\xi^{[k]}=(\xi^1_{[t_1]},...,\xi^k_{[t_k]}),& w^{[k]}=(w^1_{[t_1]},...,w^k_{[t_k]})\nonumber
\end{eqnarray}
the collection of all random vectors at iteration $k$ and the collection of all such random vectors from iteration $1$ to $k$ respectively, and note that given $(\xi^{[k]},w^{[k]})$, $z^{k+1}$ is a deterministic vector, we then have the following bound:
\begin{eqnarray*}
&& \mathbb{E}\left[(z^{k+1}-z^*)^\top\left(\hat{F}^{k+1}_{\rho}(z^{k+1})-F(z^{k+1})\right)\right]\\
&=& \mathbb{E}\left[\mathbb{E}\left[(z^{k+1}-z^*)^\top\left(\hat{F}^{k+1}_{\rho}(z^{k+1})-F(z^{k+1})\right)\left|\xi^{[k]},w^{[k]}\right.\right]\right]\\
&\overset{\eqref{batch-mean}}{=}& \mathbb{E}\left[(z^{k+1}-z^*)^\top\left(F_{\rho}(z^{k+1})-F(z^{k+1})\right)\right]\\
&\ge& -\mathbb{E}\left[\|z^{k+1}-z^*\|\left\|F_{\rho}(z^{k+1})-F(z^{k+1})\right\|\right]\\
&\overset{\eqref{smooth-grad-bd}}{\ge}& -\mathbb{E}\left[\|z^{k+1}-z^*\|\cdot\frac{L\sqrt{\rho_x^2n^2+\rho_y^2m^2}}{2}\right]\\
&\ge& -\mathbb{E}\left[\frac{\mu}{2}\|z^{k+1}-z^*\|^2+\frac{L^2(\rho_xn^2+\rho_y^2m^2)}{8\mu}\right],\nonumber\\
&=& -\frac{\mu}{2}d_{k+1}-\frac{L^2(\rho_x^2n^2+\rho_y^2m^2)}{8\mu}.\nonumber
\end{eqnarray*}

Substituting the above bound into \eqref{szo-exm-bd-1}, the desired result follows.

\subsection{Proof of Proposition \ref{prop:szo-exm-sample}}
\label{proof:szo-exm-sample}
Let us start from the potential function inequality \eqref{szo-exm-potential-relation} from iteration $K$. With $\theta=\frac{1}{8}$, let us also denote $C_1=\left(1+\frac{1}{8\kappa}\right)^{-1}=\left(1-\frac{1}{8\kappa+1}\right)$, then $\rho_x=\frac{\sqrt{C_1^K}}{\sqrt{2}n\kappa},\rho_y=\frac{\sqrt{C_1^K}}{\sqrt{2}m\kappa}$. Note that the $C_1$ defined here is only for this proof, not to be confused with that used in the proof in Appendix \ref{proof:prop-szo-exp-sample}. Then we have
\begin{eqnarray}
V_K&\le& C_1V_{K-1}+16C_1\tau^2\sigma^2\left(\frac{1}{t_{K-1}}+\frac{1}{t_{K-2}}\right)+L^2\left(4\tau^2+\frac{\alpha}{8\mu}\right)\cdot C_1\cdot \frac{C_1^K}{\kappa^2}\nonumber\\
&\le& C_1^KV_0+48\tau^2\tilde{\sigma}^2\sum\limits_{k=0}^{K-1}\frac{C_1^{K-k}}{t_k}+L^2\left(4\tau^2+\frac{\alpha}{8\mu}\right)\cdot\sum\limits_{k=1}^{K}C_1^{k} \cdot \frac{C_1^K}{\kappa^2}. \nonumber
\end{eqnarray}

In the second inequality, we take $t_k=K\cdot C_1^{-k}$ and note that $\frac{1}{t_{k-1}}=\frac{1}{C_1t_k}\le\frac{2}{t_k}$. Then we have $\sum\limits_{k=0}^{K-1}\frac{C_1^{K-k}}{t_k}=C_1^K$. In addition,
\[
\sum\limits_{k=1}^KC_1^k=\frac{C_1(1-C_1^K)}{1-C_1}\le \frac{C_1}{1-C_1}=8\kappa.
\]
Therefore, we have
\begin{eqnarray}
V_K\le C_1^KV_0+48\tau^2\tilde{\sigma}^2\cdot C_1^K+L^2\left(32\tau^2+\frac{\alpha}{\mu}\right)\frac{C_1^K}{\kappa}.\label{szo-exm-potential-bd-final-1}
\end{eqnarray}
Now, let us lower bound $V_k$ by observing:
\begin{eqnarray*}
&&\tau(z^k-z^*)^\top\left(\hat{F}^{k-1}_{\rho}(z^{k-1})-\hat{F}^k_{\rho}(z^{k})\right)\\
&\ge& -\tau \|z^k-z^*\|\|\hat{F}^{k-1}_{\rho}(z^{k-1})-\hat{F}^k_{\rho}(z^{k})\|\\
&\ge& -\frac{1}{4}\|z^k-z^*\|^2-\tau^2\|\hat{F}^{k-1}_{\rho}(z^{k-1})-\hat{F}^k_{\rho}(z^{k})\|^2\\
&\overset{\eqref{Lip-szo-grad}}{\ge}& -\frac{1}{4}\|z^k-z^*\|^2-\tau^2\left(2L^2\|z^{k-1}-z^k\|^2+4(\varepsilon_{z^{k-1}}+\varepsilon_{z^k})^2+4L^2(\rho_x^2n^2+\rho_y^2m^2)\right),
\end{eqnarray*}
Then we have
\begin{eqnarray}
V_K&\ge& \frac{1}{4}d_K-16\tau^2\tilde{\sigma}^2\left(\frac{1}{t_K}+\frac{1}{t_{K-1}}\right)-\frac{4\tau^2L^2C_1^K}{\kappa^2}\nonumber\\
&\ge& \frac{1}{4}d_K-48\tau^2\tilde{\sigma}^2C_1^K-\frac{4\tau^2L^2C_1^K}{\kappa^2}.\nonumber
\end{eqnarray}

Combining with \eqref{szo-exm-potential-bd-final-1}, we have
\[
d_K\le C_1^K\cdot\left(4V_0+384\tau^2\tilde{\sigma}^2+L^2\left(\frac{32\tau^2}{\kappa}+\frac{4\tau^2}{\kappa^2}+\frac{\alpha}{\mu\kappa}\right)\right)=C_1^K\mathcal{O}(1).
\]
It follows immediately that for $K=\mathcal{O}\left(\kappa\ln\left(\frac{1}{\epsilon}\right)\right)$ we have $d_K\le\epsilon$. With a more precise expression $K=\ln_{C_1^{-1}}\left(\frac{1}{\epsilon}\right)$, the sample complexity can be estimated:
\begin{eqnarray}
\sum\limits_{k=0}^{K-1}t_k&=&K\cdot \frac{1-C_1^{-K}}{C_1(1-C_1^{-1})} = \frac{K(C_1^{-K}-1)}{1-C_1}\le \frac{\kappa}{\epsilon(1-C_1)}\ln\left(\frac{1}{\epsilon}\right)=\mathcal{O}\left(\frac{\kappa^2}{\epsilon}\ln\left(\frac{1}{\epsilon}\right)\right).\nonumber
\end{eqnarray}

\section{Proof of the uniform sublinear convergence of the stochastic extra-point method}
\label{appendix:sublinear-conv}

In order to establish a uniform sublinear convergence, we have to consider parameters that are diminishing with iteration number $k$. Let us return to the one-iteration relation \eqref{sto-exp-relation} and consider the following choice of parameters:
\[
(\alpha^{(k)},\beta^{(k)},\gamma^{(k)},\eta^{(k)},\tau^{(k)})=\left(\frac{2}{(k+2)\mu},\frac{\alpha^2\mu^2}{128},\frac{\alpha^2\mu^2}{128},\frac{2}{(k+2)\mu},\frac{\alpha^2\mu}{128\kappa}\right),
\]
where we omit the superscript $(k)$ of $\alpha$ on the RHS for notation simplicity. We shall note that here $\alpha=\alpha^{(k)}$ which is dependent on iteration $k$ and follow the same simplification for other parameters throughout the rest of the proof in this appendix unless noted otherwise.

By using the fact $2\alpha\le \frac{1}{\mu}+\alpha^2\mu$, we have:
\begin{eqnarray}
&&\left(2\alpha^2L^2+2|\gamma-\beta|+2\gamma+2\alpha\mu-1\right)\mathbb{E}\left[\|z^{k+0.5}-z^k\|^2\right]\nonumber\\
&\le& \left(2\alpha^2L^2+2\gamma+\alpha^2\mu\right)\mathbb{E}\left[\|z^{k+0.5}-z^k\|^2\right] \nonumber\\
&\le& \alpha^2\left(2L^2+\frac{\mu}{64}+\mu\right)\mathbb{E}\left[\|z^{k+0.5}-z^k\|^2\right] \nonumber\\
&\le& \alpha^2\left(2L^2+\frac{\mu^2}{64}+\mu^2\right)D^2,\nonumber
\end{eqnarray}
where in the last inequality we use the boundedness of the feasible set.

Therefore, we could rewrite \eqref{sto-exp-relation} into:
\begin{eqnarray}
&&(1-\tau L)d_{k+1}\nonumber\\
&\le& \left(1+4\gamma+4\tau L-\alpha\mu\right)d_k+\left(2\gamma+4\tau L\right)d_{k-1}\nonumber\\
&&+\alpha^2\left(2L^2+\frac{\mu^2}{64}+\mu^2\right)D^2+8\left(\alpha^2+\frac{\tau}{L}\right)\sigma^2+2\alpha\delta D\nonumber\\
&=& \left(1+4\gamma+4\tau L-\alpha\mu\right)d_k+\left(2\gamma+4\tau L\right)d_{k-1}\nonumber\\
&&+ \frac{4}{(k+2)^2}\cdot\underbrace{\left(2\kappa^2D^2+\frac{D^2}{64}+D^2+\frac{8\sigma^2}{\mu^2}+\frac{\sigma^2}{128L^2}\right)}_{G}+\frac{4\delta D}{(k+2)\mu}.\nonumber
\end{eqnarray}
Substituting the parameters with their respective values in the rest of the terms:
\begin{eqnarray}
&&\left(1-\frac{1}{32(k+2)}\right)d_{k+1}\nonumber\\
&\le& \left(1-\frac{1}{32(k+2)^2}\right)d_{k+1}\nonumber\\
&\le&\left(1+\frac{1}{8(k+2)^2}+\frac{1}{8(k+2)^2}-\frac{2}{k+2}\right)d_k\nonumber\\
&&+\left(\frac{1}{16(k+2)^2}+\frac{1}{8(k+2)^2}\right)d_{k-1}+\frac{4G}{(k+2)^2}+\frac{4\delta D}{(k+2)\mu}\nonumber\\
&\le& \left(1-\frac{7}{4(k+2)}\right)d_k+\frac{3}{16(k+2)}d_{k-1}+\frac{4G}{(k+2)^2}+\frac{4\delta D}{(k+2)\mu}.\nonumber
\end{eqnarray}

Dividing both sides by $1-\frac{1}{32(k+2)}$, and noting that $\left(1-\frac{1}{32(k+2)}\right)^{-1}\le\frac{32}{31}$, it follows that
\begin{eqnarray}
d_{k+1}&\le& \frac{1-\frac{7}{4(k+2)}}{1-\frac{1}{32(k+2)}}d_k+\frac{6}{31(k+2)}d_{k-1}+\frac{128G}{31(k+2)^2}+\frac{128\delta D}{31(k+2)\mu}\nonumber\\
&=& \left(1-\frac{\frac{55}{32(k+2)}}{1-\frac{1}{32(k+2)}}\right)d_k+\frac{6}{31(k+2)}d_{k-1}+\frac{128G}{31(k+2)^2}+\frac{128\delta D}{31(k+2)\mu}\nonumber\\
&\le& \left(1-\frac{55}{32(k+2)}\right)d_k+\frac{7}{32(k+2)}d_{k-1}+\frac{128G}{31(k+2)^2}+\frac{128\delta D}{31(k+2)\mu}.\nonumber
\end{eqnarray}
From the above one-iteration inequality, we shall claim the following inequality
\[
d_k\le \frac{Q}{k+2}+\frac{256\delta D}{93\mu},\quad \forall k\ge0,
\]
where
\[
Q=\max\left\{\frac{133G}{9},2\|z^0-z^*\|^2\right\},
\]
and we shall prove the inequality by induction. For $k=0$, the inequality holds trivially
\[
d_0=\|z^0-z^*\|^2\le\frac{Q}{2}.
\]
Assuming the inequality holds for all index $1,...,k$, we then have
\begin{eqnarray}
d_{k+1}&\le& \left(1-\frac{55}{32(k+2)}\right)d_k+\frac{7}{32(k+2)}d_{k-1}+\frac{128G}{31(k+2)^2}+\frac{128\delta D}{31(k+2)\mu}\nonumber\\
&\le& \left(1-\frac{55}{32(k+2)}\right)\left(\frac{Q}{k+2}+\frac{256\delta D}{93\mu}\right)+\frac{7}{32(k+2)}\left(\frac{Q}{k+1}+\frac{256\delta D}{93\mu}\right)\nonumber\\
&&+\frac{128G}{31(k+2)^2}+\frac{128\delta D}{31(k+2)\mu}\nonumber\\
&\le&\left(1-\frac{55}{32(k+2)}\right)\cdot\frac{Q}{k+2}+\frac{7}{32(k+2)}\cdot\frac{Q}{k+1}+\frac{128G}{31(k+2)^2}\nonumber\\
&&+\frac{256\delta D}{93\mu}-\frac{440\delta D}{93(k+2)\mu}+\frac{56\delta D}{93(k+2)\mu}+\frac{128\delta D}{31(k+2)\mu}\nonumber\\
&=& \frac{Q}{k+2}-\frac{55Q}{32(k+2)^2}+\frac{7}{32(k+2)}\cdot\frac{Q}{k+1}+\frac{128G}{31(k+2)^2}+\frac{256\delta D}{93\mu}\nonumber\\
&\le& \frac{Q}{k+3}+\frac{Q}{(k+2)^2}-\frac{55Q}{32(k+2)^2}+\frac{7}{32(k+2)}\cdot\frac{2Q}{k+2}+\frac{133G}{32(k+2)^2}+\frac{256\delta D}{93\mu}\nonumber\\
&=& \frac{Q}{k+3}-\frac{9Q}{32(k+2)^2}+\frac{133G}{32(k+2)^2}+\frac{256\delta D}{93\mu}\nonumber\\
&=&\frac{Q}{k+3}+\frac{256\delta D}{93\mu}.\nonumber
\end{eqnarray}
Note that in the last inequality we used the identities $\frac{1}{k+1}\le\frac{2}{k+2}$ and $\frac{128}{31}\le\frac{133}{32}$. This completes the proof for the uniform $\mathcal{O}\left(\frac{1}{k}\right)$ convergence rate.

\end{appendices}

\end{document}